\documentclass[11pt]{article}

\usepackage{amsmath,amssymb, amsthm}
\usepackage{fullpage}

\usepackage{amsfonts}
\usepackage{layout}
\usepackage{amsmath}
\usepackage{amsthm}
\usepackage{amssymb}
\usepackage{url}
\usepackage{color}
\usepackage{graphicx}
\usepackage{verbatim}

\newcommand{\eqdef}{\stackrel{\text{def}}{=}}

\newcommand{\R}{\mathbb{R}}
\newcommand{\Prob}{\mathbf{Prob}}
\newcommand{\E}{\mathbf{E}}

\theoremstyle{plain}
\newtheorem{theorem}{Theorem}

\newtheorem{lemma}[theorem]{Lemma}

\theoremstyle{definition}

\makeatletter
\newcommand*{\rom}[1]{\expandafter\@slowromancap\romannumeral #1@}
\makeatother

\usepackage{algorithm}
\usepackage[noend]{algpseudocode}
\usepackage{subcaption}

\title{Primal Method for ERM  with  Flexible Mini-batching Schemes \\and Non-convex Losses\footnote{The authors acknowledge support from the EPSRC Grant EP/K02325X/1, Accelerated Coordinate Descent Meth- ods for Big Data Optimization.}}

\author{Dominik Csiba \thanks{School of Mathematics, The University of Edinburgh, United Kingdom (e-mail: cdominik@gmail.com)}\quad  \qquad \quad Peter Richt\'arik \thanks{School of Mathematics, The University of Edinburgh, United Kingdom (e-mail: peter.richtarik@ed.ac.uk)}}

\begin{document}
\maketitle

\begin{abstract}
In this work we develop a new algorithm for regularized empirical risk minimization. Our method extends recent techniques of Shalev-Shwartz [02/2015], which enable a dual-free analysis of SDCA, to arbitrary mini-batching schemes.  Moreover, our method is able to better utilize the information in the data defining the ERM problem. For convex loss functions, our complexity results  match those of QUARTZ, which is  a primal-dual method also allowing for arbitrary mini-batching schemes. The advantage of a dual-free analysis comes from the fact that it guarantees convergence even for non-convex loss functions, as long as the average loss is convex. We illustrate through experiments the utility of being able to design arbitrary mini-batching schemes.
\end{abstract}

\section{Introduction}

Empirical risk minimization (ERM) is a very successful and immensely popular paradigm in machine learning, used to train a variety of prediction and classification models. Given  examples $A_1,\dots,A_n\in \R^{d\times m}$, loss functions $\phi_1,\dots,\phi_n:\R^m\to \R$ and a regularization parameter $\lambda>0$, the L2-regularized ERM problem is an optimization problem of the form
\begin{equation} \label{eq:gen_prob}
\min_{w \in \R^d} \left[ P(w) := \frac{1}{n} \sum_{i=1}^n \phi_i(A_i^\top w) + \frac{\lambda}{2} \|w\|^2 \right]
\end{equation}

Throughout the paper we shall assume that for each $i$, the loss function $\phi_i$ is $l_i$-smooth with $l_i>0$. That is, for all $ x,y \in \R^m$ and all $i\in [n]:=\{1,2,\dots,n\}$, we have 
\begin{equation} \label{eq:phiAiw_lsmoothness}
\|\nabla\phi_i(x) - \nabla\phi_i(y)\| \leq l_i\|x - y\|.
\end{equation}

Further, let $L_1,\dots, L_n>0$ be constants for which the inequality 
\begin{equation} \label{eq:phiAiw_smoothness}
\|\nabla\phi_i(A_i^\top w) - \nabla\phi_i(A_i^\top z)\| \leq L_i\|w - z\|
\end{equation}
holds for all  $w,z \in \R^d$ and all $i$ and let $L := \max_i L_i$. Note that we can always bound $L_i \leq l_i \|A_i\|$. However, $L_i$ can be better (smaller) than $l_i \|A_i\|$.

\subsection{Background} In the last few years, a lot of research effort was put into  designing new efficient algorithms for solving this problem (and some of its modifications). The frenzy of activity was motivated by the realization that SGD  \cite{robbins1951}, not so long ago considered the state-of-the-art method for ERM, was far from being optimal, and that new ideas can lead to algorithms which are far superior to SGD in both theory and practice. The methods that belong to this category include  SAG \cite{SAG}, SDCA \cite{SDCA}, SVRG \cite{SVRG}, S2GD \cite{S2GD}, mS2GD \cite{mS2GD}, SAGA \cite{SAGA},  S2CD \cite{S2CD}, QUARTZ \cite{Quartz}, ASDCA \cite{ASDCA}, prox-SDCA \cite{PROXSDCA}, IPROX-SDCA \cite{IProx-SDCA}, A-PROX-SDCA \cite{A-PROX-SDCA}, AdaSDCA \cite{AdaSDCA}, SDNA \cite{SDNA}. Methods analyzed for arbitrary mini-batching schemes include 
NSync \cite{NSync},  ALPHA \cite{ALPHA} and QUARTZ \cite{Quartz}.


In order to find an $\epsilon$-solution in expectation, state of the art (non-accelerated) methods for solving \eqref{eq:gen_prob} only need  \[O((n+\kappa)\log(1/\epsilon))\] steps, where each step involves the computation of the gradient $\nabla \phi_i(A_i^\top w)$ for some randomly selected example $i$. The quantity $\kappa$ is the condition number. Typically one has $\kappa = \frac{\max_i l_i \|A_i\|^2}{\lambda}$ for methods picking $i$ uniformly at random, and $\kappa = \frac{\sum_i l_i \|A_i\|^2}{n\lambda}$ for methods picking $i$ using  a carefully designed data-dependent importance sampling. Computation of such a  gradient  typically involves work which is equivalent to reading the example $A_i$, that is, $O(nnz(A_i))\leq dm$ arithmetic operations.

\subsection{Contributions} In this work we develop a new algorithm for the L2-regularized ERM problem \eqref{eq:gen_prob}.  Our method extends a technique recently introduced by  Shalev-Shwartz  \cite{DualFreeSDCA}, which enables a dual-free analysis of SDCA, to {\em arbitrary mini-batching schemes.}   That is, our method works at each iteration with a random subset of examples, chosen in an i.i.d. fashion from an arbitrary distribution. Such flexible schemes are useful for various reasons, including i) the development of distributed or robust variants of the method, ii) design of importance sampling for improving the complexity rate, iii) design of a sampling which is aimed at obtaining efficiencies elsewhere, such us utilizing NUMA (non-uniform memory access) architectures, and  iv) streamlining and speeding up the processing of each mini-batch by means of assigning to each processor approximately even workload so as to reduce idle time (we do experiments with the latter setup).  

In comparison with \cite{DualFreeSDCA},  our method is able to better {\em utilize the information in the data examples} $A_1,\dots,A_n$, leading to a better {\em data-dependent} bound.  For convex loss functions, our complexity results  match those of QUARTZ \cite{Quartz} in terms of the rate (the logarithmic factors differ). QUARTZ  is  a primal-dual method also allowing for arbitrary mini-batching schemes. However, while  \cite{Quartz} only characterize the decay of expected risk, we also give bounds for the {\em sequence of iterates}. In particular, we show that for convex loss functions, our method enjoys the rate (Theorem~\ref{thm:phiAiw_convex})
\[\max_i \left(\frac{1}{p_i}  +  \frac{l_i v_i}{\lambda p_i n}\right) \log\left(\frac{(L + \lambda)E^{(0)}}{\lambda\epsilon}\right),\] where $p_i$ is the probability that coordinate $i$ is updated in an iteration,   $v_1,\dots,v_n>0$ are certain ``stepsize'' parameters of the method associated with the sampling and data (see \eqref{eq:phiAiw_ESO}), and $E^{(0)}$ is a constant depending on the starting point.  For instance, in the special case picking a single example at a time uniformly at random, we have $p_i=1/n$ and $v_i=\|A_i\|^2$, whereby we obtain one of the $O(n+\kappa)\log(1/\epsilon)$ rates mentioned above. The other rate can be recovered using importance sampling. 

The advantage of a dual-free analysis comes from the fact that it guarantees convergence even for {\em non-convex} loss functions, as long as the average loss is convex. This is a step toward understanding non-convex models. In particular, we show that for non-convex loss functions, our method enjoys the rate (Theorem~\ref{thm:phiAiw_nonconvex})
\[  \max_i \left(  \frac{1}{p_i} + \frac{L_i^2 v_i}{\lambda^2 p_i n}\right) \log\left(\frac{(L + \lambda)D^{(0)}}{\lambda\epsilon}\right),\] 
where $D^{(0)}$ is a constant depending on the starting point.

Finally, we illustrate through experiments with ``chunking''---a simple load balancing technique---the utility of being able to design arbitrary mini-batching schemes.


\section{Algorithm} \label{sec:phiAiw}

We shall now describe the method (Algorithm~\ref{alg:phiAiw}). 

\begin{algorithm}
\caption{dfSDCA: Dual-Free SDCA with Arbitrary Sampling} \label{alg:phiAiw}
\begin{algorithmic}
\State \textbf{Parameters:} Sampling $\hat{S}$, stepsize $\theta$
\State \textbf{Initialization} $\alpha_1^{(0)}, \dots, \alpha_n^{(0)} \in \R^m$, set $w^{(0)} = \frac{1}{\lambda n}\sum_{i=1}^n A_i\alpha_i^{(0)},~p_i = \mathbf{Prob}(i \in \hat{S})$
\For{$t \geq 1$}
\State Sample a set $S_t$ according to $\hat{S}$
\For{$i \in S_t$}
\State $\alpha_i^{(t)} = \alpha_i^{(t-1)} - \theta p_i^{-1}(\nabla\phi_i(A_i^\top w^{(t-1)}) + \alpha_i^{(t-1)})$
\EndFor
\State $w^{(t)} = w^{(t-1)} - \sum_{i\in S_t}\theta (n \lambda p_i)^{-1} A_i(\nabla\phi_i(A_i^\top w^{(t-1)}) + \alpha_i^{(t-1)})$ 
\EndFor
\end{algorithmic}
\end{algorithm}

The method encodes a family of algorithms, depending on the choice of the sampling $\hat{S}$, which encodes a particular {\em mini-batching scheme}.  Formally, a sampling $\hat{S}$ is a set-valued random variable with values being the subsets of $[n]$, i.e., subsets of examples. In this paper, we use the terms ``mini-batching scheme'' and ``sampling'' interchangeably. A sampling is defined by the collection of probabilities $\Prob(S)$  assigned to every subset $S\subseteq [n]$ of the examples.


The method maintains $n$ vectors $\alpha_i\in \R^m$ and a vector $w\in \R^d$. At the beginning of  step $t$, we have $\alpha_i^{(t-1)}$ for all $i$ and $w^{(t-1)}$ computed and stored in memory. We then pick a random subset $S_t$ of the examples, according to the mini-batching scheme, and update variables $\alpha_i$ for $i \in S_t$, based on the computation of the gradients $\nabla \phi_i(A_i^\top w^{(t-1)})$ for $i\in S_t$. This is followed by an update of the vector $w$, which is performed so as to maintain the relation
\begin{equation}\label{eq:rel}
w^{(t)} = \frac{1}{\lambda n}\sum_i A_i \alpha_i^{(t)}.\end{equation} This relation is maintained for the following reason. If  $w^*$ is the optimal solution to \eqref{eq:gen_prob}, then
\begin{equation}
0 = \nabla P(w^*) = \frac{1}{n} \sum_{i=1}^n A_i \nabla \phi_i (A_i^\top w^*) + \lambda w^* ,\label{eq:phiAiw_optimality}
\end{equation}
and hence
$ w^* = \frac{1}{\lambda n} \sum_{i=1}^n A_i \alpha_i^*$,
where $\alpha_i^* := -\nabla\phi_i(A_i^\top w^*)$. So, if we believe that the variables $\alpha_i$ converge to $-\nabla\phi_i(A_i^\top w^*)$, it indeed does make sense to maintain \eqref{eq:rel}. Why should we believe this? This is where the specific update of the ``dual variables'' $\alpha_i$ comes from: $\alpha_i$ is set a convex combination of its previous value and our best estimate so far of $ -\nabla\phi_i(A_i^\top w^*)$, namely, $-\nabla \phi_i(A_i^\top w^{(t-1)})$. Indeed, the update can be written as
\[\alpha_i^{(t)} = (1-\theta p_i^{-1})\alpha_i^{(t-1)} + \theta p_i^{-1}(-\nabla \phi_i(A_i^\top w^{(t-1)})).\] Why does {\em this} make sense? Because we believe that $w^{(t-1)}$ converges to $w^*$.  Admittedly, this  reasoning is somewhat ``circular''. However, a better word to describe this reasoning would be: ``iterative''.



\section{Main Results}

Let $p_i := \mathbb{P} (i \in \hat{S})$. We assume the knowledge of parameters $v_1, \dots, v_n>0$ for which
\begin{equation} \label{eq:phiAiw_ESO}
 \E \left[ \left\| \sum_{i \in \hat{S}} A_i h_i \right\|^2 \right] \leq \sum_{i=1}^n p_i v_i \|h_i\|^2.
\end{equation}
Tight and easily computable formulas for such parameters can be found in \cite{ESO}.  For instance, whenever $\Prob(|\hat{S}|\leq \tau)=1$, inequality \eqref{eq:phiAiw_ESO} holds with $v_i = \tau \|A_i\|^2$.

To simplify the exposure, we will write \begin{equation} \label{def:phiAitw_pot}
 \qquad B^{(t)} \eqdef \|w^{(t)} - w^*\|^2, \qquad C_i^{(t)} \eqdef \| \alpha_i^{(t)} - \alpha_i^*\|^2, \qquad i=1,2,\dots,n.
\end{equation}

\subsection{Non-convex loss functions}

Our result will be expressed in terms of the decay of the  potential
$
D^{(t)} \eqdef \frac{\lambda}{2} B^{(t)} + \frac{\lambda}{2n} \sum_{i=1}^n \frac{1}{L_i^2}C_i^{(t)}$,
where $B_i^{(t)}$ and $C^{(t)}$ are defined in \eqref{def:phiAitw_pot}.
\begin{theorem} \label{thm:phiAiw_nonconvex}
Assume that  the average loss function, $\frac{1}{n}\sum_{i=1}^n \phi_i$, is convex. If  \eqref{eq:phiAiw_smoothness} holds and we let \begin{equation} \label{eq:phiAiw_eta}
\theta \leq  \min_i \frac{ p_i n \lambda^2}{L_i^2 v_i + n \lambda^2},
\end{equation} 
then  the for $t\geq 0$ the potential $D^{(t)}$ decays exponentially to zero as
\begin{equation}\E \left[D^{(t)}\right]
 \leq e^{-\theta t} D^{(0)} .
\end{equation}
Moreover, if we set $\theta$ equal to the upper bound in \eqref{eq:phiAiw_eta},  then 
\[
T \geq  \max_i \left(  \frac{1}{p_i} + \frac{L_i^2 v_i}{\lambda^2 p_i n}\right) \log\left(\frac{(L + \lambda)D^{(0)}}{\lambda\epsilon}\right) \quad \Rightarrow \quad  \E[P(w^{(T)}) - P(w^*)] \leq \epsilon.\]
\end{theorem}

\subsection{Convex loss functions}

 Our result will be expressed in terms of the decay of the  potential
$E^{(t)} \eqdef \frac{\lambda}{2}B^{(t)} + \frac{1}{2n}\sum_{i=1}^n \frac{1}{l_i}C_i^{(t)}$, where $B_i^{(t)}$ and $C^{(t)}$ are defined in \eqref{def:phiAitw_pot}. 

\begin{theorem} \label{thm:phiAiw_convex}
Assume that all loss functions $\{\phi_i\}$ are convex and satisfy \eqref{eq:phiAiw_lsmoothness}. If we run Algorithm~\ref{alg:phiAiw} with parameter $\theta$ satisfying the inequality
\begin{equation} \label{eq:phiAiw_eta_convex}
\theta \leq \min_i \frac{p_i n \lambda}{l_i v_i + n \lambda},
\end{equation}
then  the for $t\geq 0$ the potential $E^{(t)}$ decays exponentially to zero as
\begin{equation}\label{eq:main}
\E \left[E^{(t)}\right]
 \leq e^{-\theta t} E^{(0)}.
\end{equation}
Moreover, if we set $\theta$ equal to the upper bound in \eqref{eq:phiAiw_eta_convex},
then 
\[T \geq  \max_i \left(\frac{1}{p_i}  + \frac{l_i v_i}{\lambda p_i n} \right) \log\left(\frac{(L+\lambda) E^{(0)}}{\lambda \epsilon}\right) \quad \Rightarrow \quad  \E[P(w^{(T)}) - P(w^*)] \leq \epsilon\]
\end{theorem}

The rate, $\theta$, precisely matches that of the QUARTZ algorithm \cite{Quartz}. Quartz is the only other method for ERM which has been analyzed for an arbitrary mini-batching scheme. Our algorithm is dual-free, and as we have seen above, allows for an analysis covering the case of non-convex loss functions.

\section{Chunking}

In this section we illustrate one use of the ability of our method to work with an arbitrary mini-batching scheme. Further examples include the ability to design distributed variants of the method \cite{Hydra}, or the use of importance/adaptive sampling to lower the number of iterations \cite{UCDC, IProx-SDCA, Quartz, AdaSDCA}. 

One marked disadvantage of standard mini-batching (``choose a subset of examples, uniformly at random'') used in the context of   parallel processing on multicore processors is the fact that in a synchronous implementation there is a loss of efficiency due to the fact that the computation time of $\nabla\phi (A_i^\top w)$ may differ through $i$. This is caused by the data examples having varying degree of sparsity. We hence introduce a new sampling which mitigates this issue.

\paragraph{Chunks:} Choose sets $G_1, \dots, G_k \subset [n]$, such that $\cup_{i=1}^k G_i = [n]$ and $G_i \cap G_j = \emptyset ~ \forall i,j$ and $\psi (i) := \sum_{j \in G_{i}} \text{nnz}(A_j)$ is similar for every $i$, i.e. $\psi (1) \approx \dots \approx \psi (k)$. Instead of sampling $\tau$ coordinates we propose a new sampling, which on each iteration $t$ samples $\tau$ sets $G_{(1)}^{(t)}, \dots, G_{(\tau )}^{(t)}$ out of $G_1, \dots, G_k$ and uses coordinates $i \in \cup_{i=1}^\tau G_{(i)}^{(t)} $ as the sampled set. We assign each core one of the sets $G_{(i)}^{(t)}$ for parallel computation. The advantage of this sampling lies in the fact, that the load of computing $\nabla \phi (A_i^\top w)$ for all  $i \in G_j$ is similar for all $j \in [k]$. Hence, using this sampling we minimize the waiting time of processors.

\paragraph{How to choose $G_1, \dots, G_k$:} We introduce the following algorithm:

\begin{algorithm}
\caption{Naive Chunks} \label{alg:naive_chunks}
\begin{algorithmic}
\State \textbf{Parameters:} vector of nnz $u$ 
\State \textbf{Initialization} $n = \text{length}(u)$; Empty vector $g$ and $s$ of length $n$; $m = \max(u)$

\State $g[1] = 1$,\quad $s[1] = u[1]$, \quad $i = 1$

\For{$t = 2:n$}
\If{$g[i] + u[t] \leq m$}
\State $g[i] = g[i] + 1$, $s[i] = s[i] + u[t]$
\Else
\State $i = i+1$, $g[i] = 1$, $s[i] = u[t]$
\EndIf

\EndFor
\end{algorithmic}
\end{algorithm}

\begin{figure}[t]      
        \centering \label{fig:chunks}
        \begin{subfigure}[b]{0.32\textwidth}
            \centering
            \includegraphics[width=\textwidth]{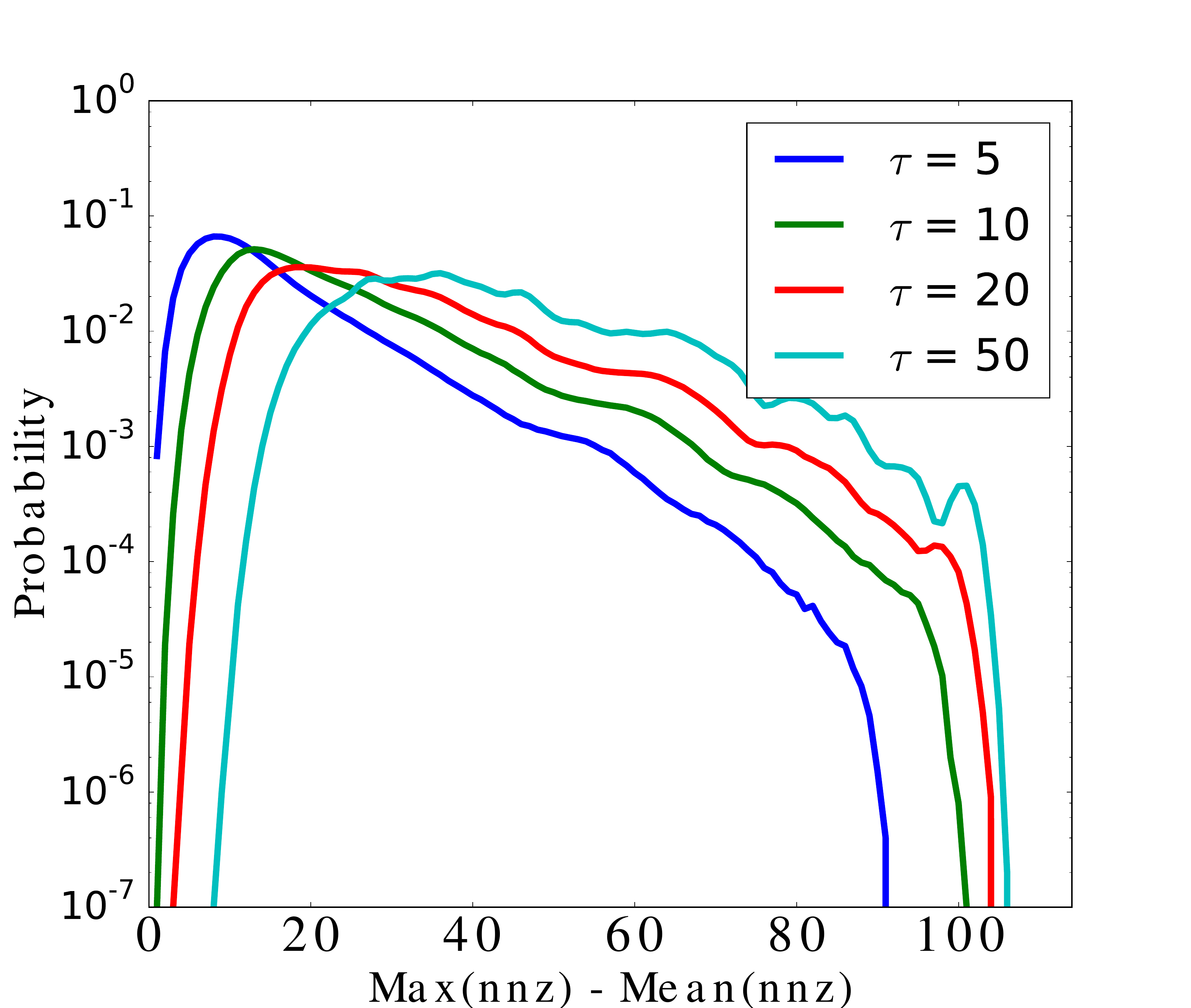}
            \caption
            {{\small w8a initially}}    
            \label{fig:chunks_initial_w8a}
        \end{subfigure}
        \begin{subfigure}[b]{0.32\textwidth}  
            \centering 
            \includegraphics[width=\textwidth]{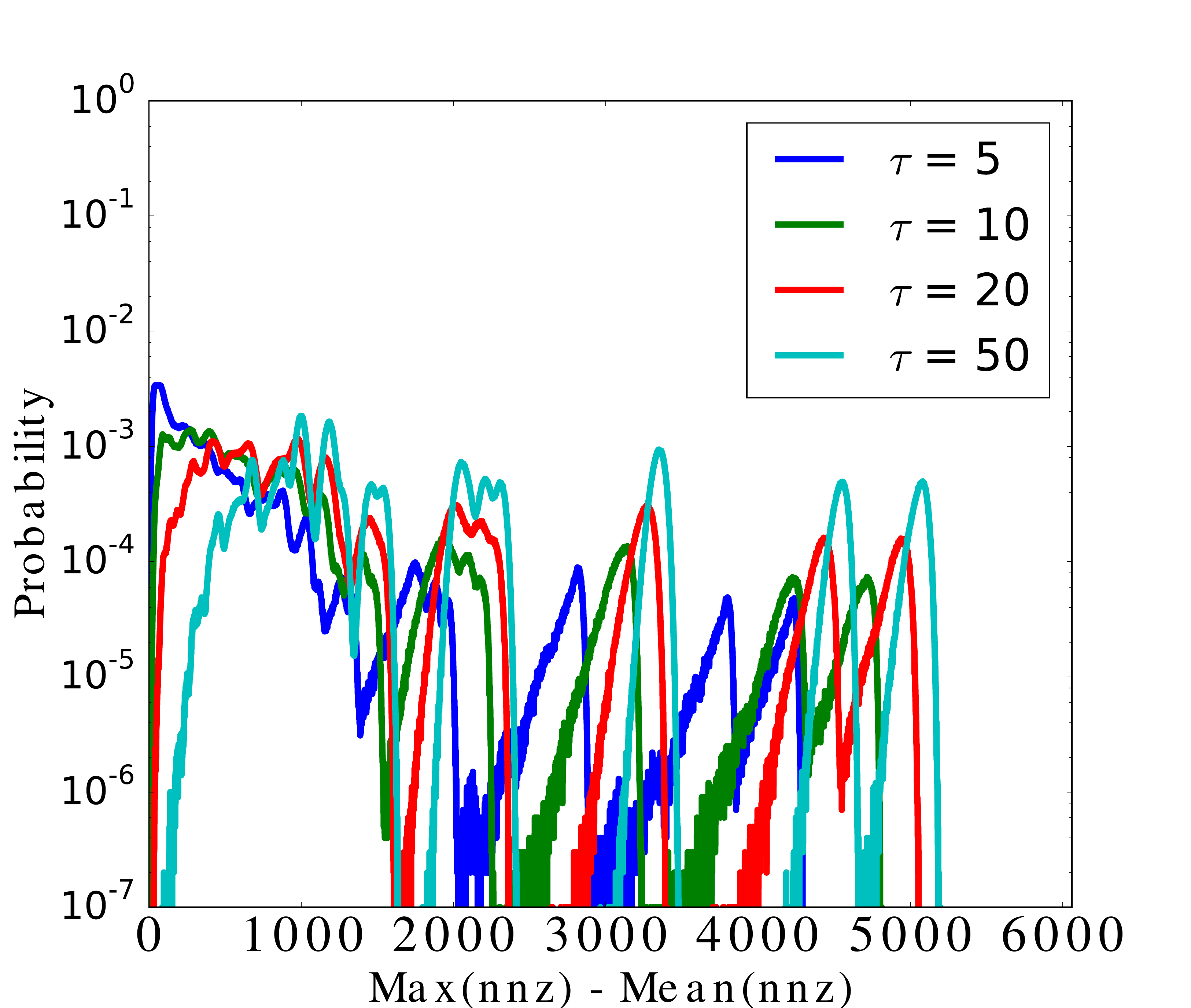}
            \caption[]%
            {{\small dorothea initially}}    
        \end{subfigure}
        \begin{subfigure}[b]{0.32\textwidth}  
            \centering 
            \includegraphics[width=\textwidth]{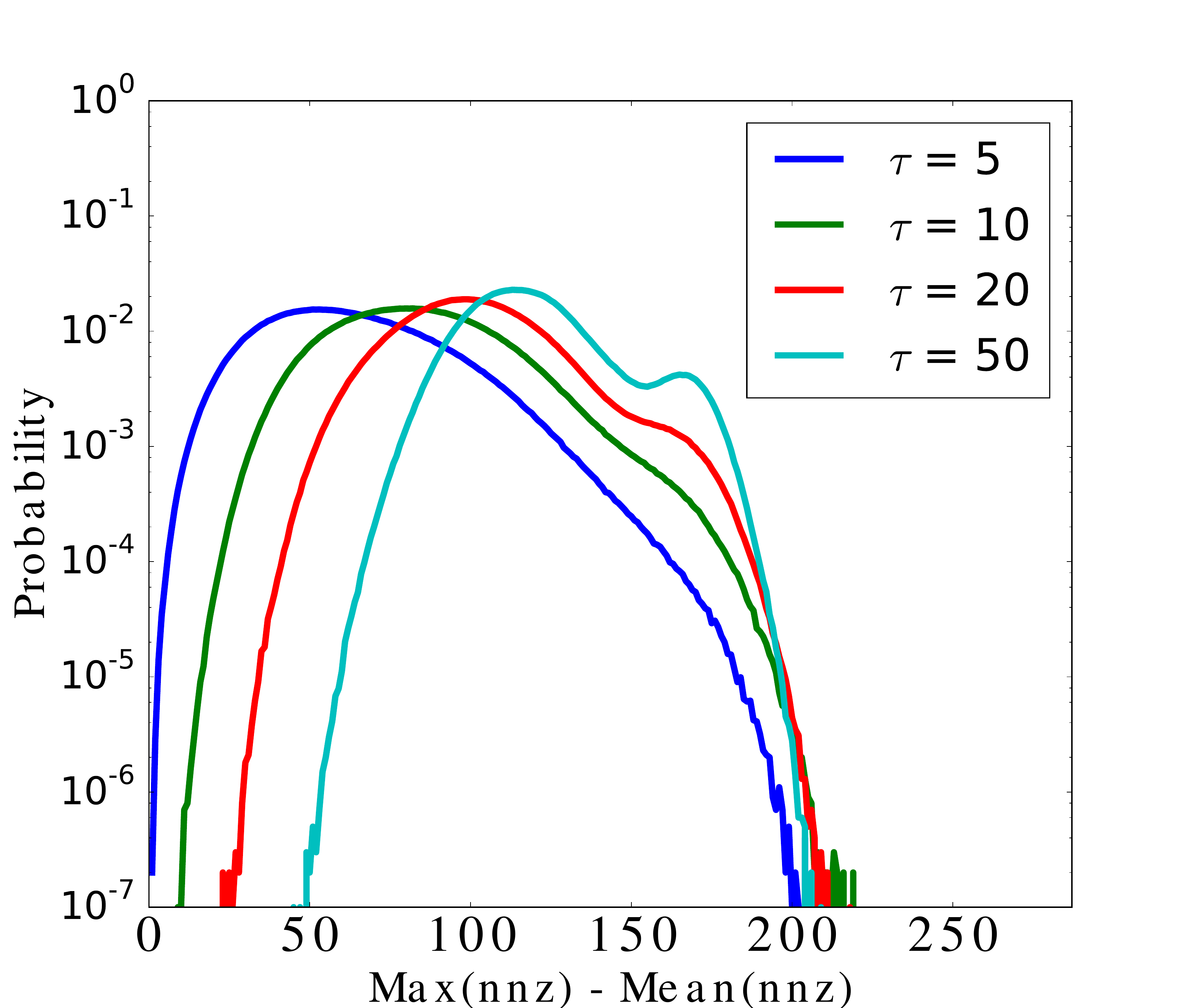}
            \caption[]%
            {{\small protein initially}}    
        \end{subfigure}
        \begin{subfigure}[b]{0.32\textwidth}   
            \centering 
            \includegraphics[width=\textwidth]{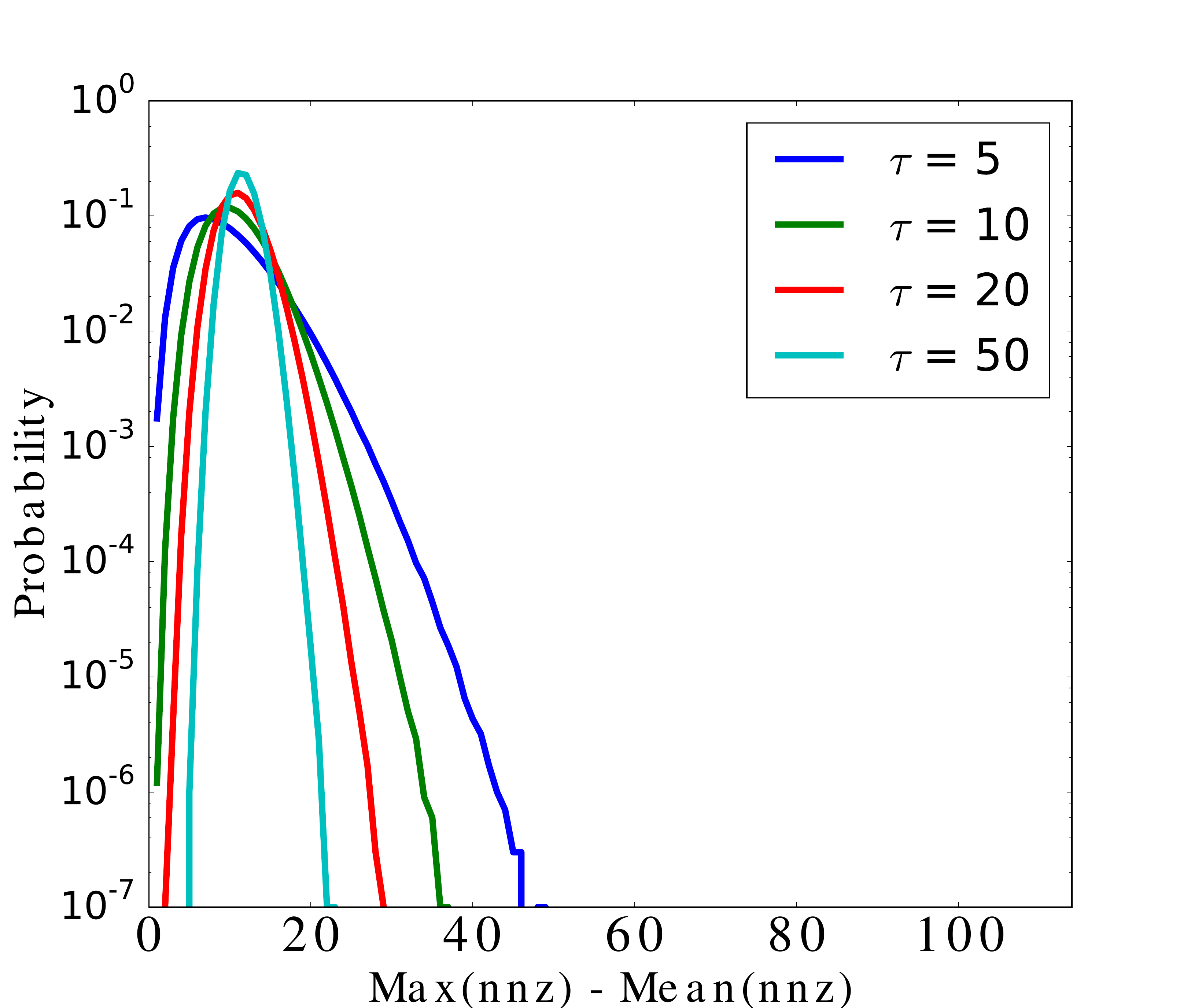}
            \caption
            {{\small w8a chunked}}    
        \end{subfigure}
        \begin{subfigure}[b]{0.32\textwidth}   
            \centering 
            \includegraphics[width=\textwidth]{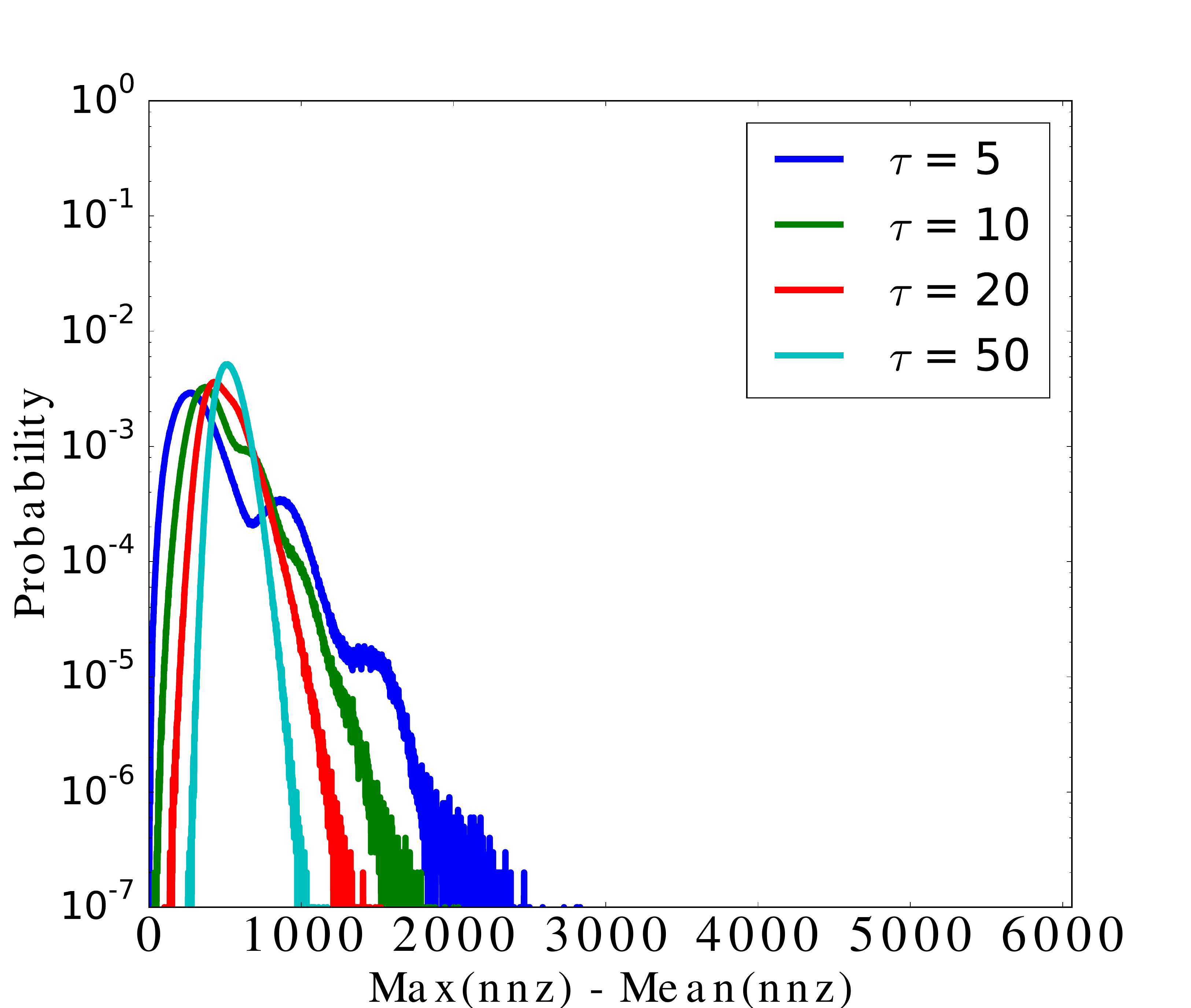}
            \caption
            {{\small dorothea chunked}}    
        \end{subfigure}
        \begin{subfigure}[b]{0.32\textwidth}  
            \centering 
            \includegraphics[width=\textwidth]{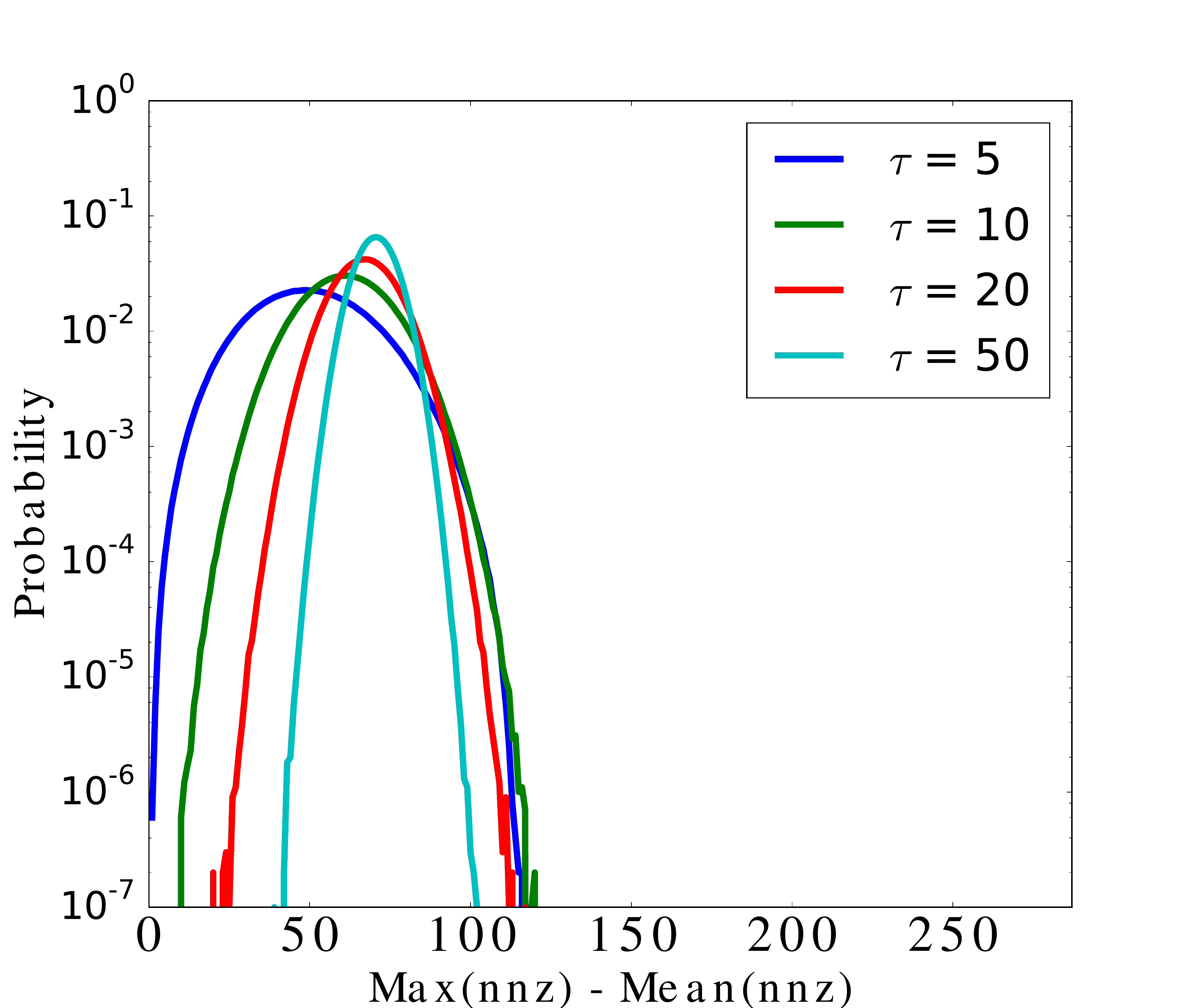}
            \caption[]%
            {{\small protein chunked}}
            \label{fig:chunks_chunked_protein}    
        \end{subfigure}
        \caption        
        {\small Distribution of the difference between the maximum number of nonzeros processed by a single core and the mean of all nonzeros processed by each core. This difference shows us, how much time is wasted per core waiting on the slowest core to finish its task, therefore smaller numbers are better. The first row corresponds to the initial distribution while the second row shows the distribution after using Algorithm~\ref{alg:naive_chunks}.} 
    \end{figure}

The algorithms returns the partition of $[n]$ into $G_1, \dots, G_k$ in a sense, that the first $g[1]$ coordinates belong to $G_1$, next $g[2]$ coordinates belong to $G_2$ and so on. The main advantage of this approach is, that it makes a preprocessing step on the dataset which takes just one pass through the data. On Figure~\ref{fig:chunks_initial_w8a} through 
Figure~\ref{fig:chunks_chunked_protein} we show the impact of Algorithm~\ref{alg:naive_chunks} on the probability of the waiting time of a single core, which we measure by the difference \[\max_{i \in S_t}\{ \text{nnz}(A_i) \} - \frac{1}{\tau}\sum_{i \in S_t} \text{nnz}(A_i)\] and \[\max_{i \in [\tau]}\{ \text{nnz}(G_{(i)}^{(t)}) \} - \frac{1}{\tau}\sum_{i=1}^\tau \text{nnz}(G_{(i)}^{(t)})\] for the initial and preprocessed dataset respectively. We can observe, that the waiting time is smaller using the preprocessing.

%

\section{Experiments}

In all our experiments we used logistic regression. We normalized the datasets so that $\max_i \|A_i\| = 1$, and fixed $\lambda = 1/n$. The datasets used for experiments are summarized in Table~\ref{tab:datasets}.

\begin{table}[!h]
\centering 
\begin{tabular}{|l|l|l|l|}
\hline
Dataset  & \#samples & \#features & sparsity \\ \hline
w8a      & 49,749    & 300        & 3.8\%    \\ \hline
dorothea & 800       & 100,000    & 0.9\%    \\ \hline
protein  & 17,766    & 358        & 29\%     \\ \hline
rcv1     & 20,242    & 47,237     & 0.2\%  	 \\ \hline
cov1	 & 581,012   & 54		  & 22\%	 \\ \hline
\end{tabular}
\caption{\footnotesize Datasets used in the experiments.}
\label{tab:datasets}
\end{table}

\begin{figure*}
        \centering
        \begin{subfigure}[b]{0.32\textwidth}
            \centering
            \includegraphics[width=\textwidth]{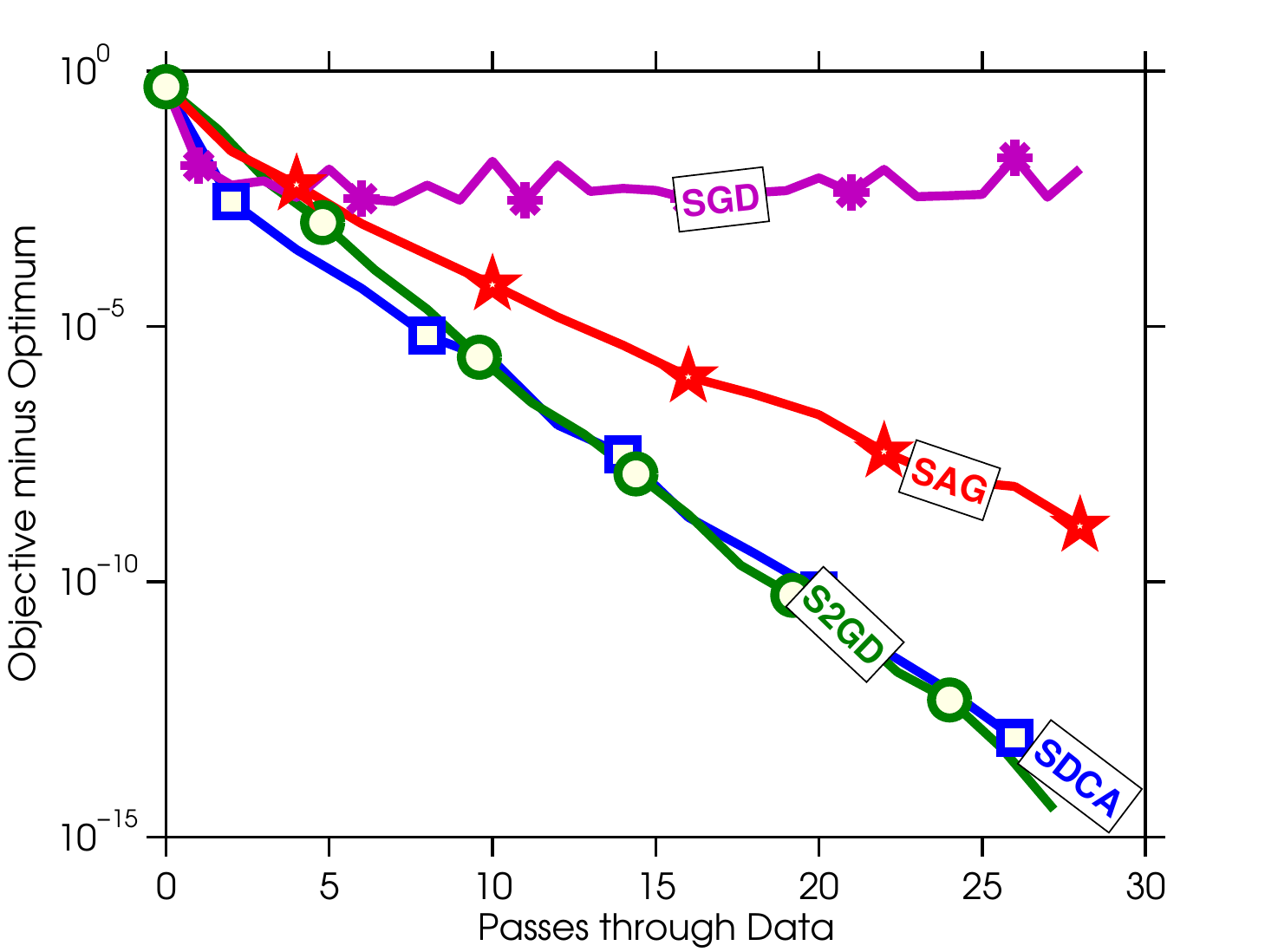}
            \caption
            {{\small rcv1, state of the art}}
            \label{fig:state_of_the_art}    
        \end{subfigure}
        \begin{subfigure}[b]{0.32\textwidth}  
            \centering 
            \includegraphics[width=\textwidth]{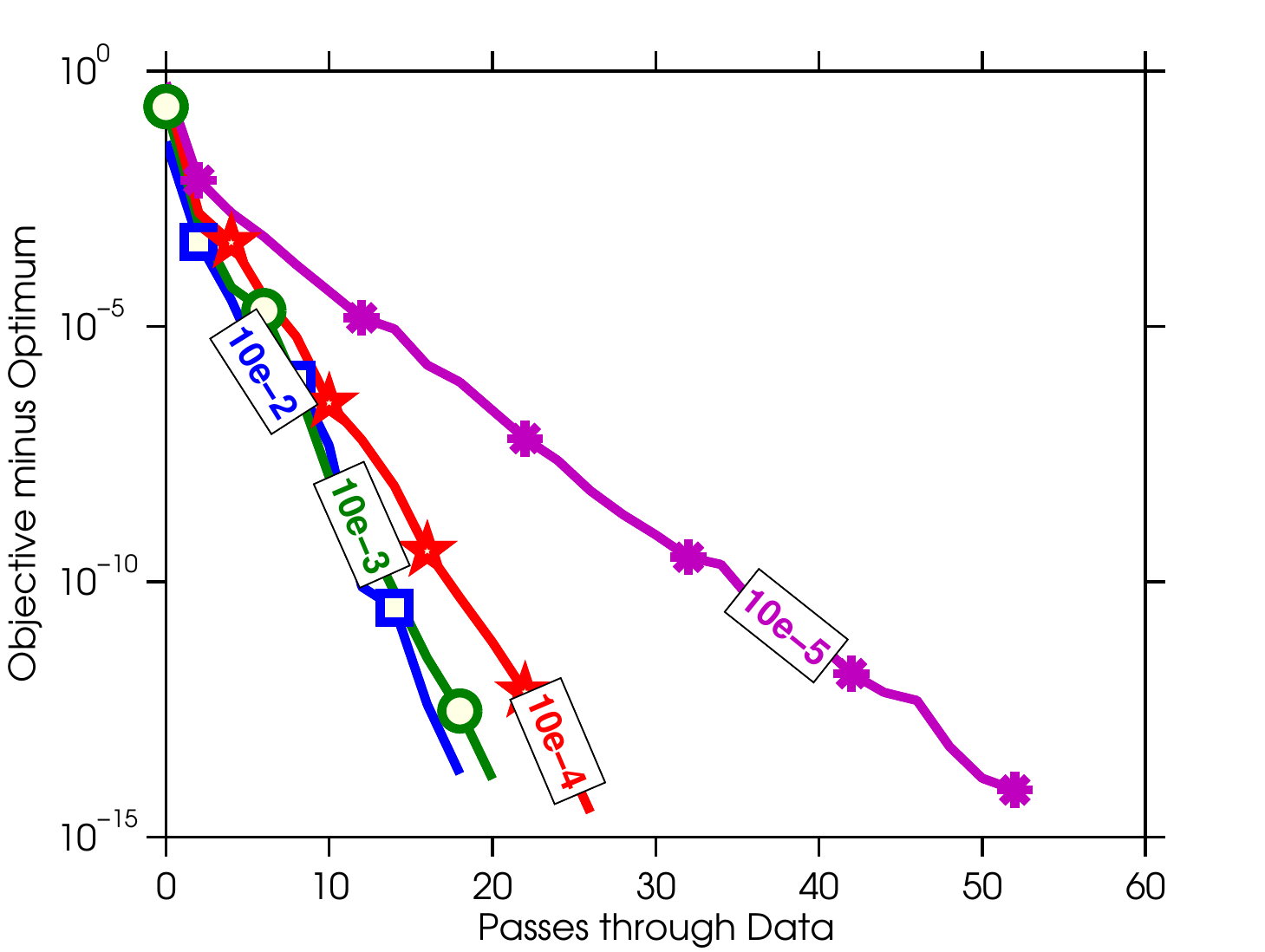}
            \caption[]%
            {{\small rcv1, different $\lambda$}}
            \label{fig:lambda}    
        \end{subfigure}
        \begin{subfigure}[b]{0.32\textwidth}  
            \centering 
            \includegraphics[width=\textwidth]{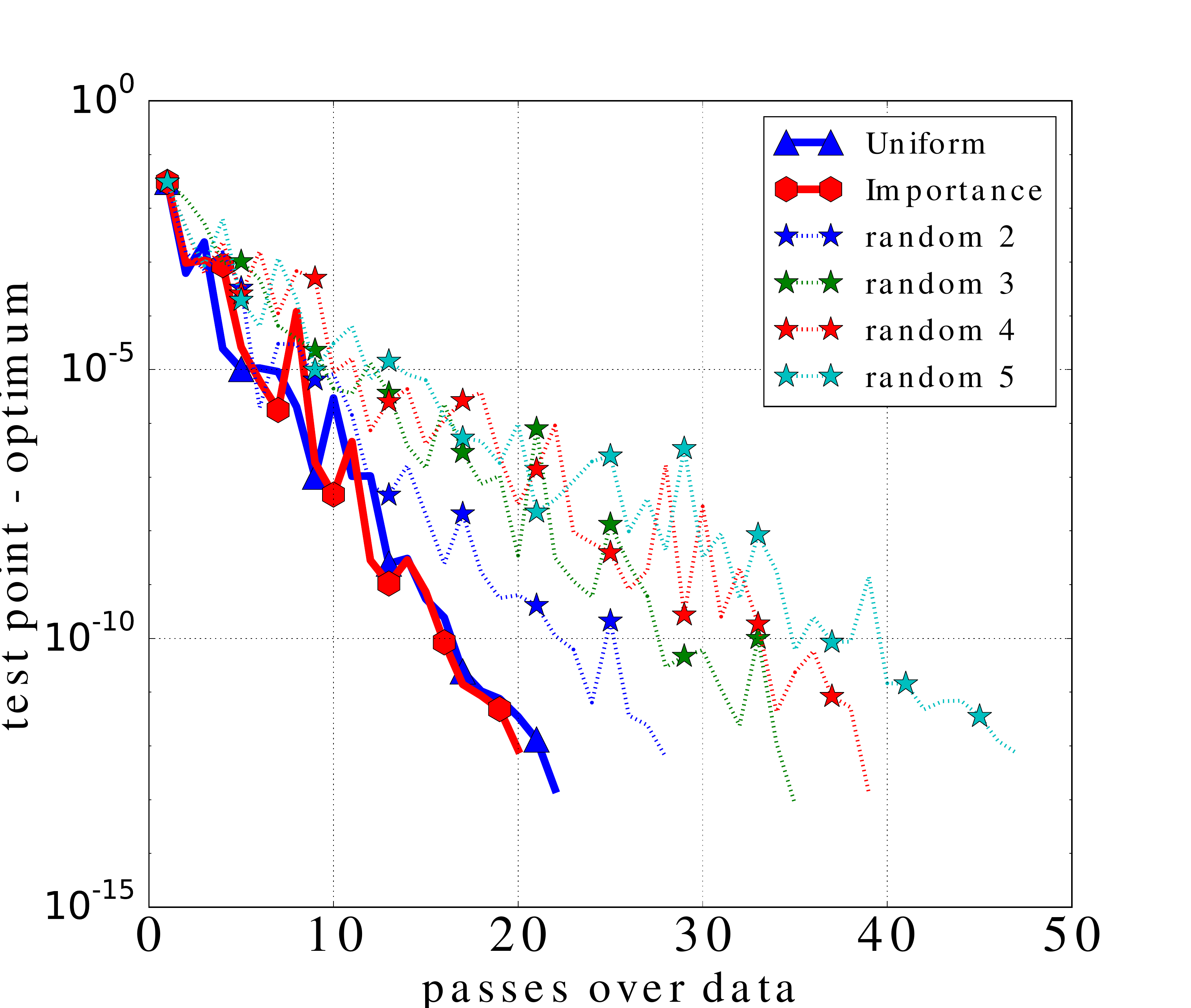}
            \caption[]%
            {{\small cov1, various samplings}}
            \label{fig:serial_samplings}    
        \end{subfigure}
 	\caption{\footnotesize LEFT: Comparison of SDCA with other state of the art methods. MIDDLE: SDCA for various values of $\lambda$. RIGHT: SDCA run with various samplings $\hat{S}$.}
    \end{figure*}

\textbf{Experiment 1.} In Figure~\ref{fig:state_of_the_art} we compared the performance of Algorithm~\ref{alg:phiAiw} with uniform serial sampling against state of the art algorithms such as SGD \cite{robbins1951}, SAG\cite{SAG} and S2GD \cite{S2GD} in number of epochs. The real running time of the algorithms was 0.46s for S2GD, 0.79s for SAG, 0.47s for SDCA and 0.58s for SGD. In Figure~\ref{fig:lambda} we show the convergence rate for different  regularization parameters $\lambda$. In Figure~\ref{fig:serial_samplings} we show convergence rates for different serial samplings: uniform, importance \cite{IProx-SDCA} and also 4 different randomly generated serial samplings. These samplings were generated in a controlled manner, such that \textit{random c} has $(\max_i~p_i)/(\min_i~p_i) < c$. All of these samplings have linear convergence as shown in the theory.

\begin{figure*}
        \centering
        \begin{subfigure}[b]{0.24\textwidth}
            \centering
            \includegraphics[width=\textwidth]{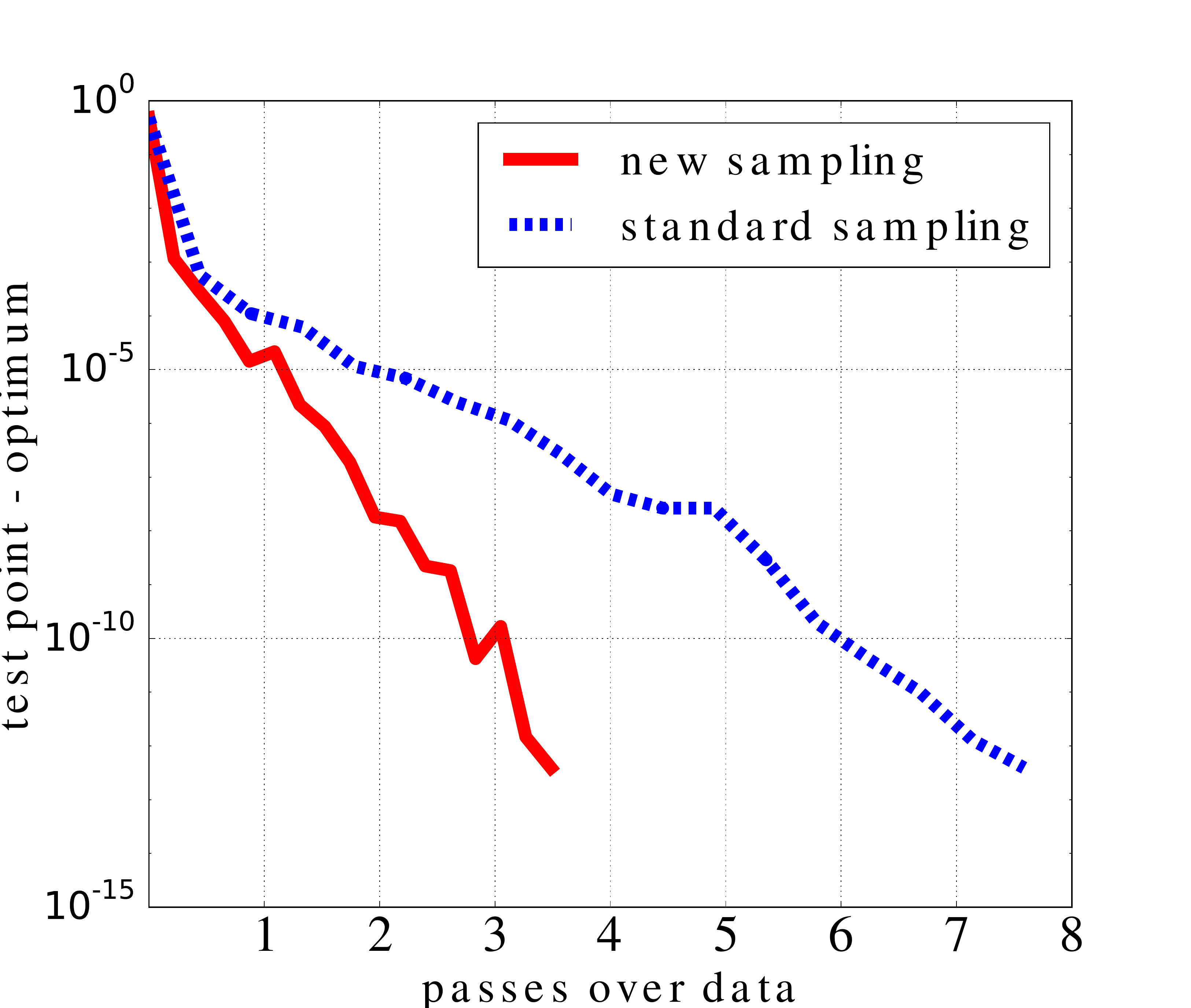}
            \caption
            {{\small w8a with $\tau= 5$}}
            \label{fig:w8a_5}    
        \end{subfigure}
        \begin{subfigure}[b]{0.24\textwidth}   
            \centering 
            \includegraphics[width=\textwidth]{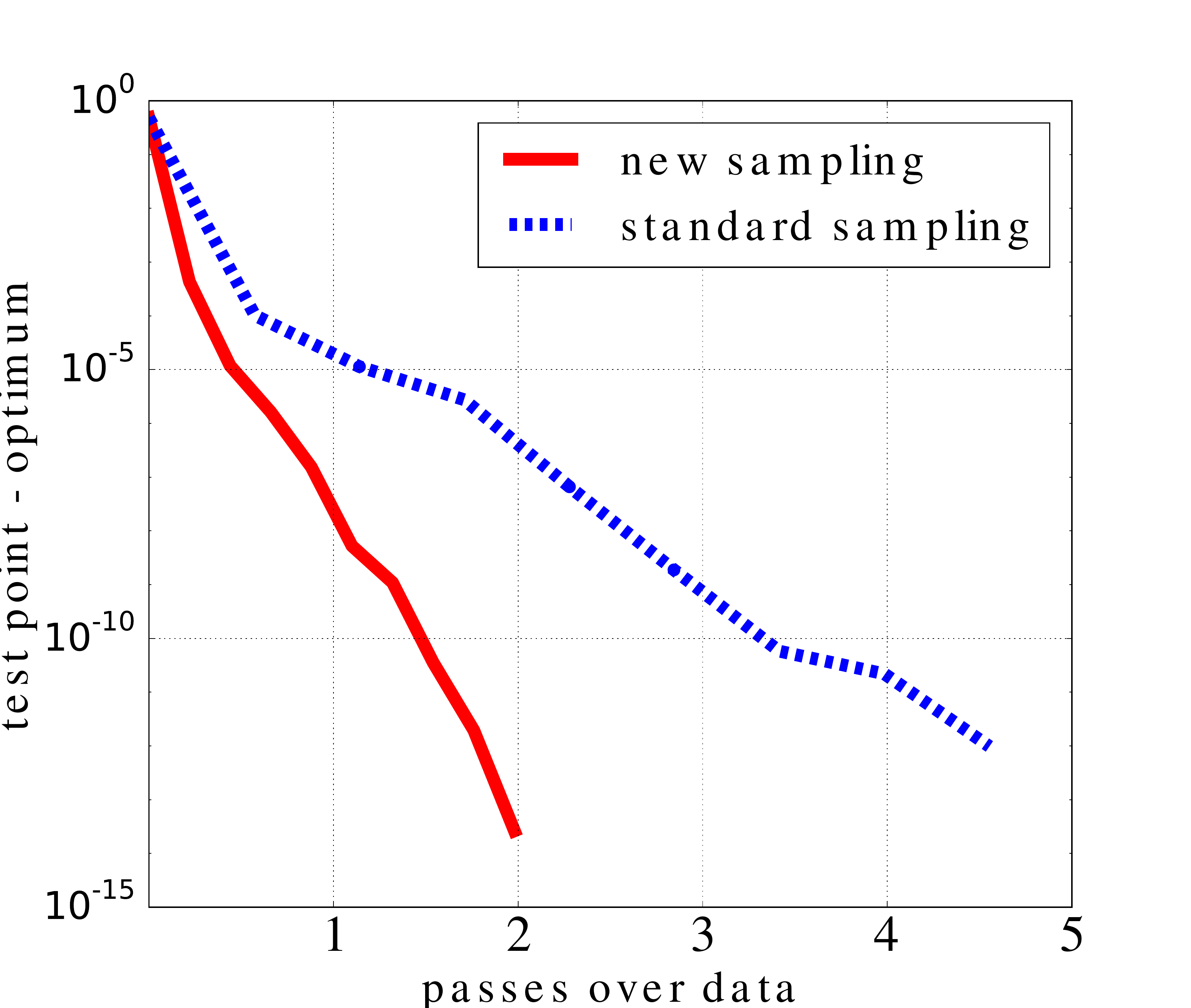}
            \caption
            {{\small w8a with $\tau= 10$}}    
        \end{subfigure}
        \begin{subfigure}[b]{0.24\textwidth}   
            \centering 
            \includegraphics[width=\textwidth]{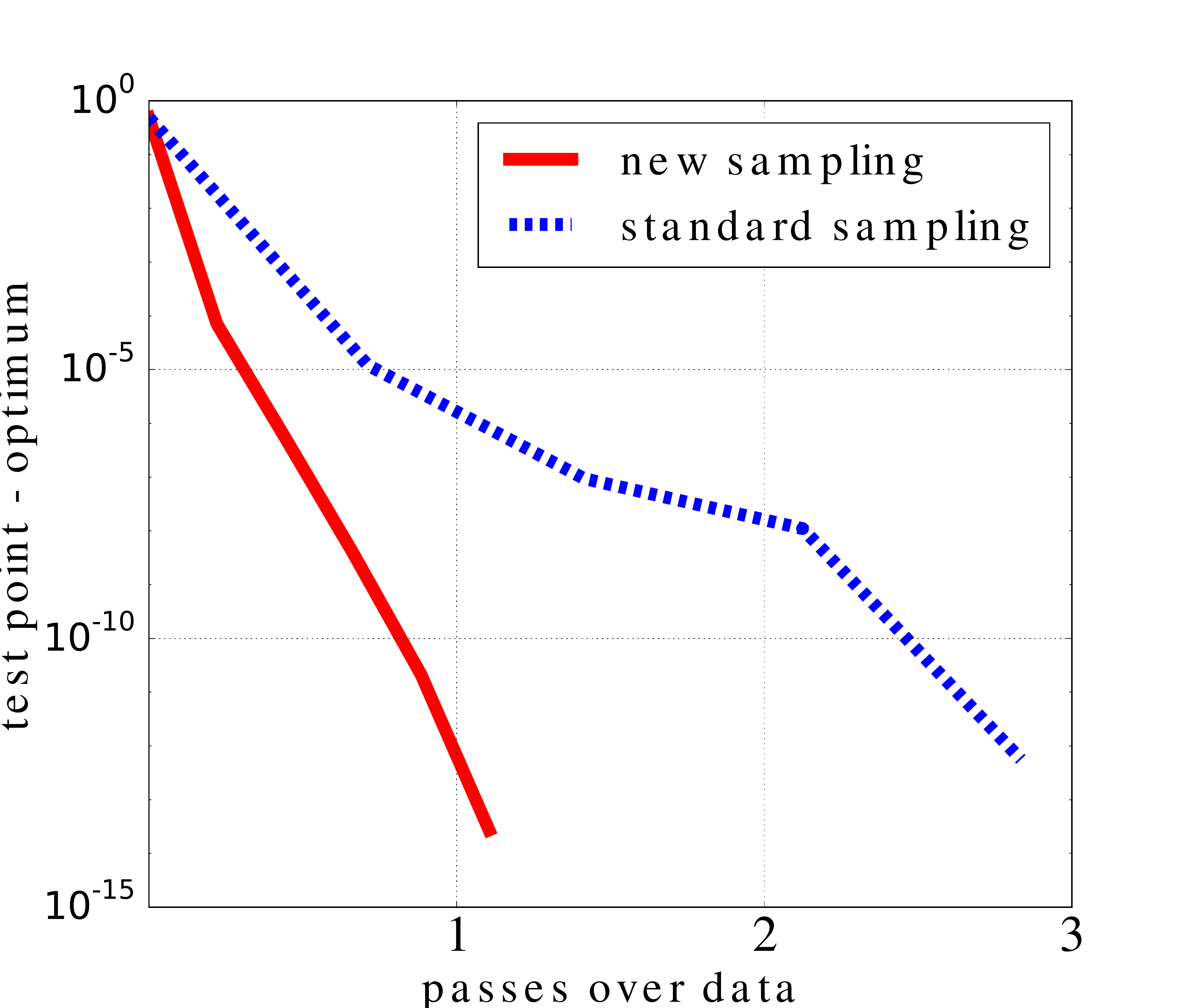}
            \caption
            {{\small w8a with $\tau= 20$}}   
        \end{subfigure}
        \begin{subfigure}[b]{0.24\textwidth}   
            \centering 
            \includegraphics[width=\textwidth]{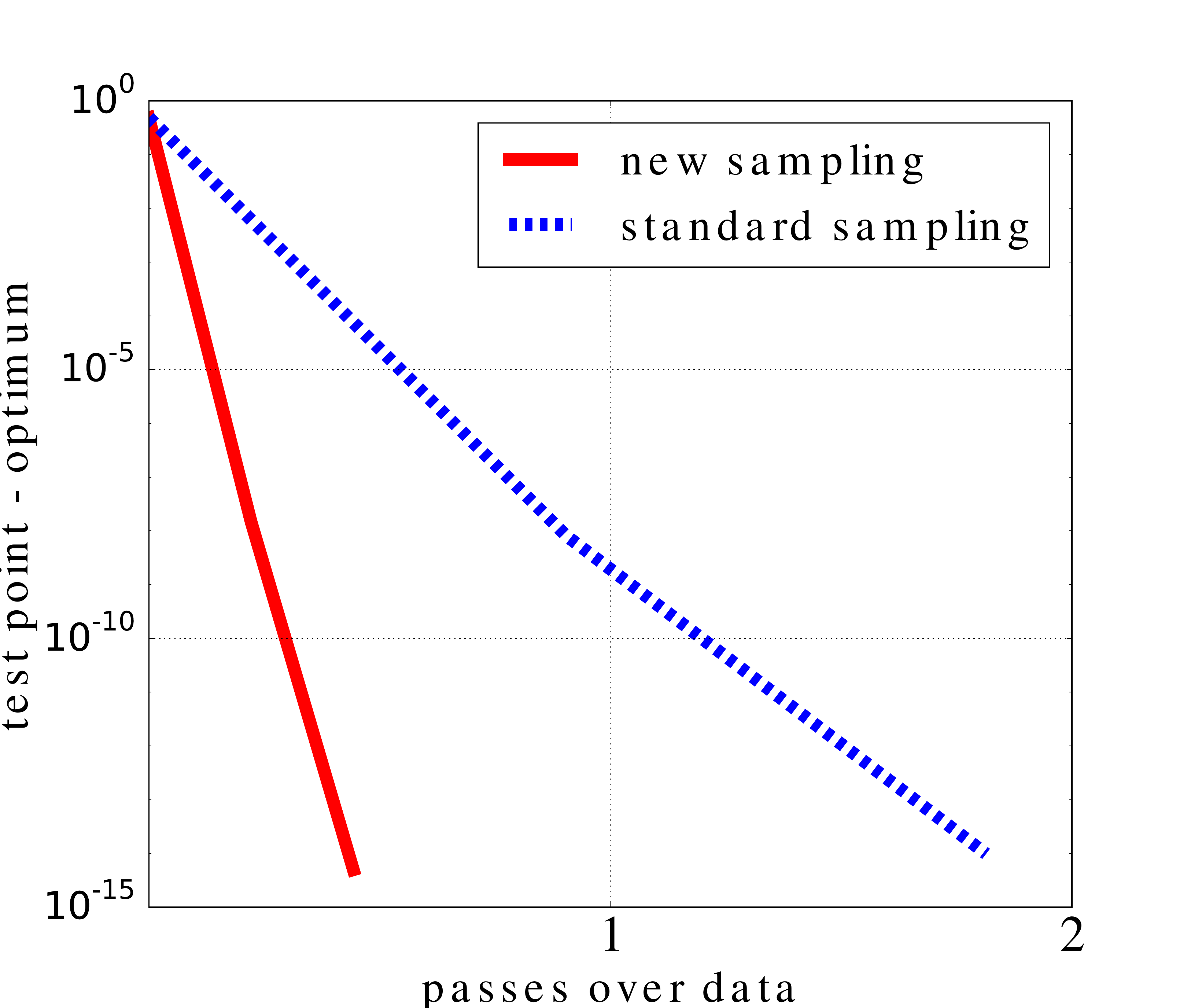}
            \caption
            {{\small w8a with $\tau= 50$}}  
        \end{subfigure}
        \begin{subfigure}[b]{0.24\textwidth}  
            \centering 
            \includegraphics[width=\textwidth]{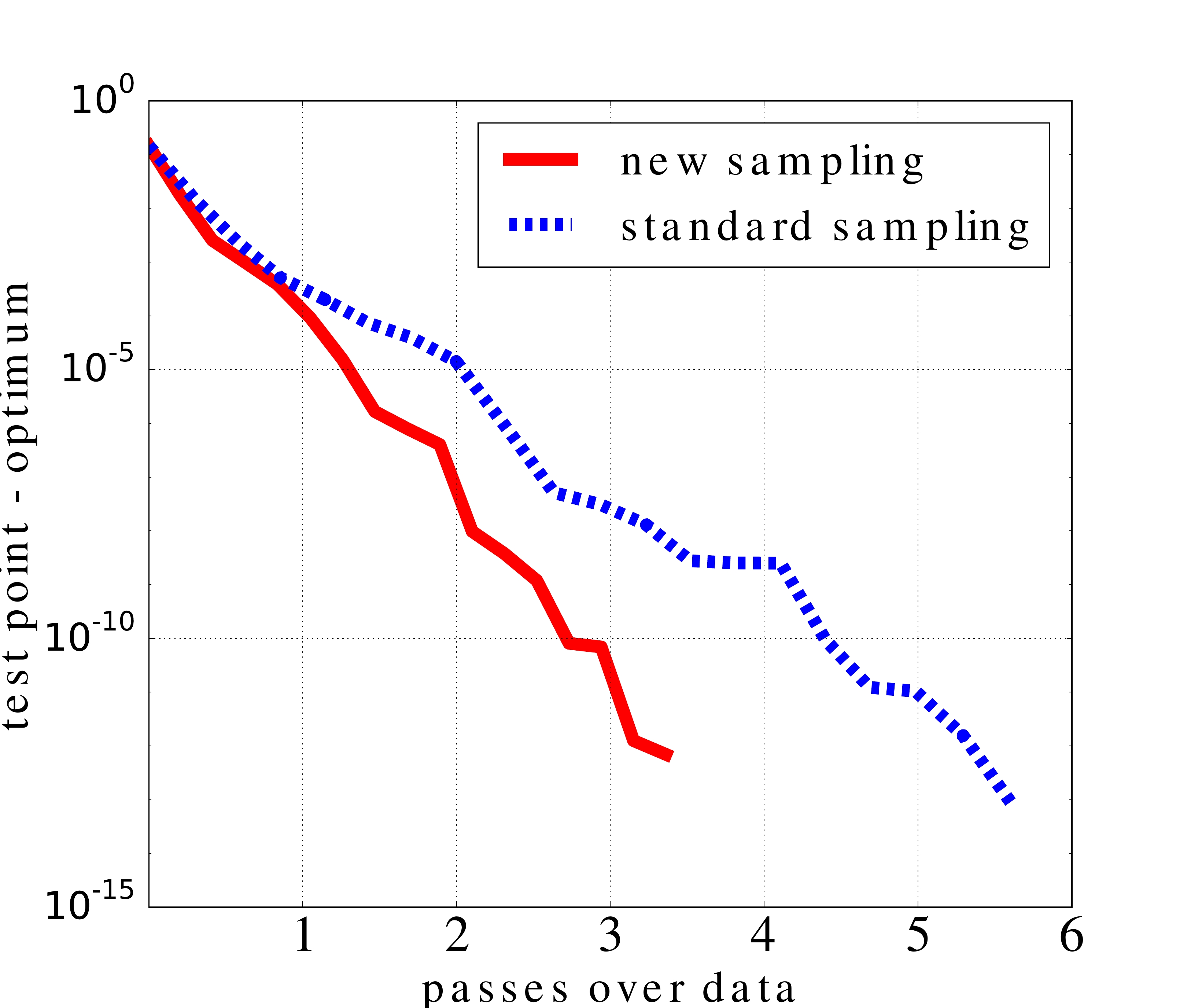}
            \caption[]%
            {{\small dorothea with $\tau=5$}}    
        \end{subfigure}                
        \begin{subfigure}[b]{0.24\textwidth}   
            \centering 
            \includegraphics[width=\textwidth]{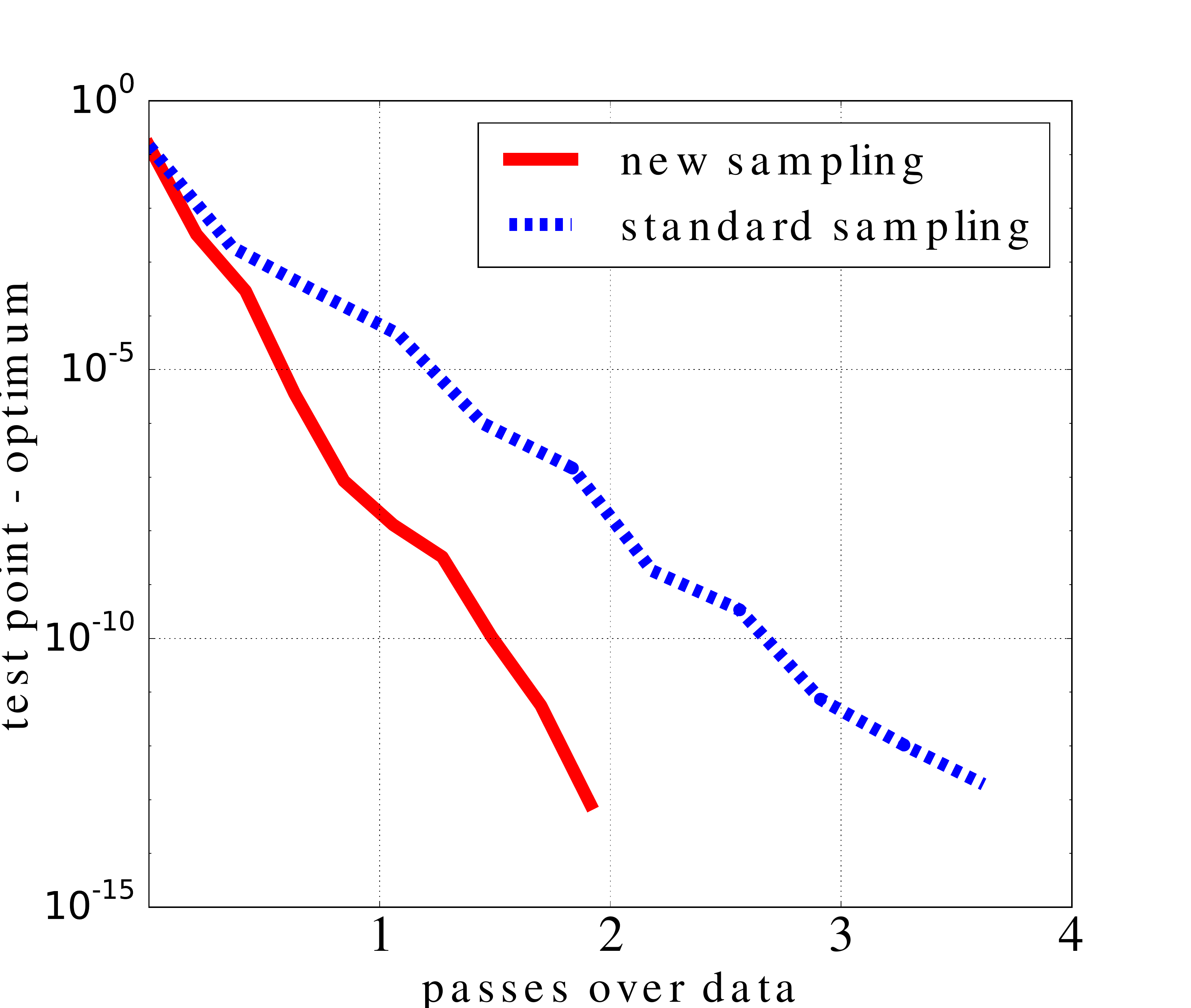}
            \caption
            {{\small dorothea with $\tau=10$}}    
        \end{subfigure}              
        \begin{subfigure}[b]{0.24\textwidth}   
            \centering 
            \includegraphics[width=\textwidth]{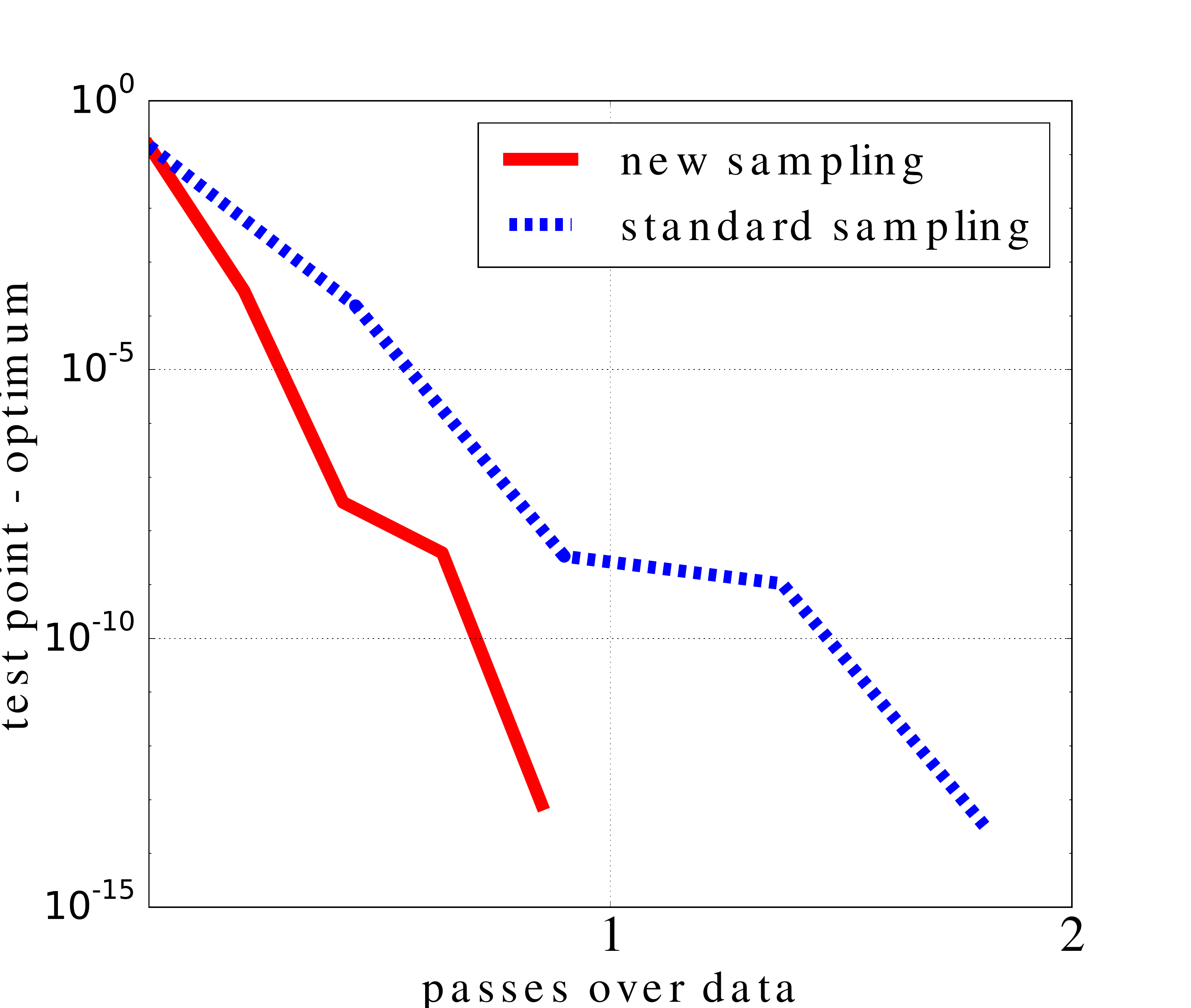}
            \caption[]%
            {{\small dorothea with $\tau=20$}}    
        \end{subfigure}      
        \begin{subfigure}[b]{0.24\textwidth}   
            \centering 
            \includegraphics[width=\textwidth]{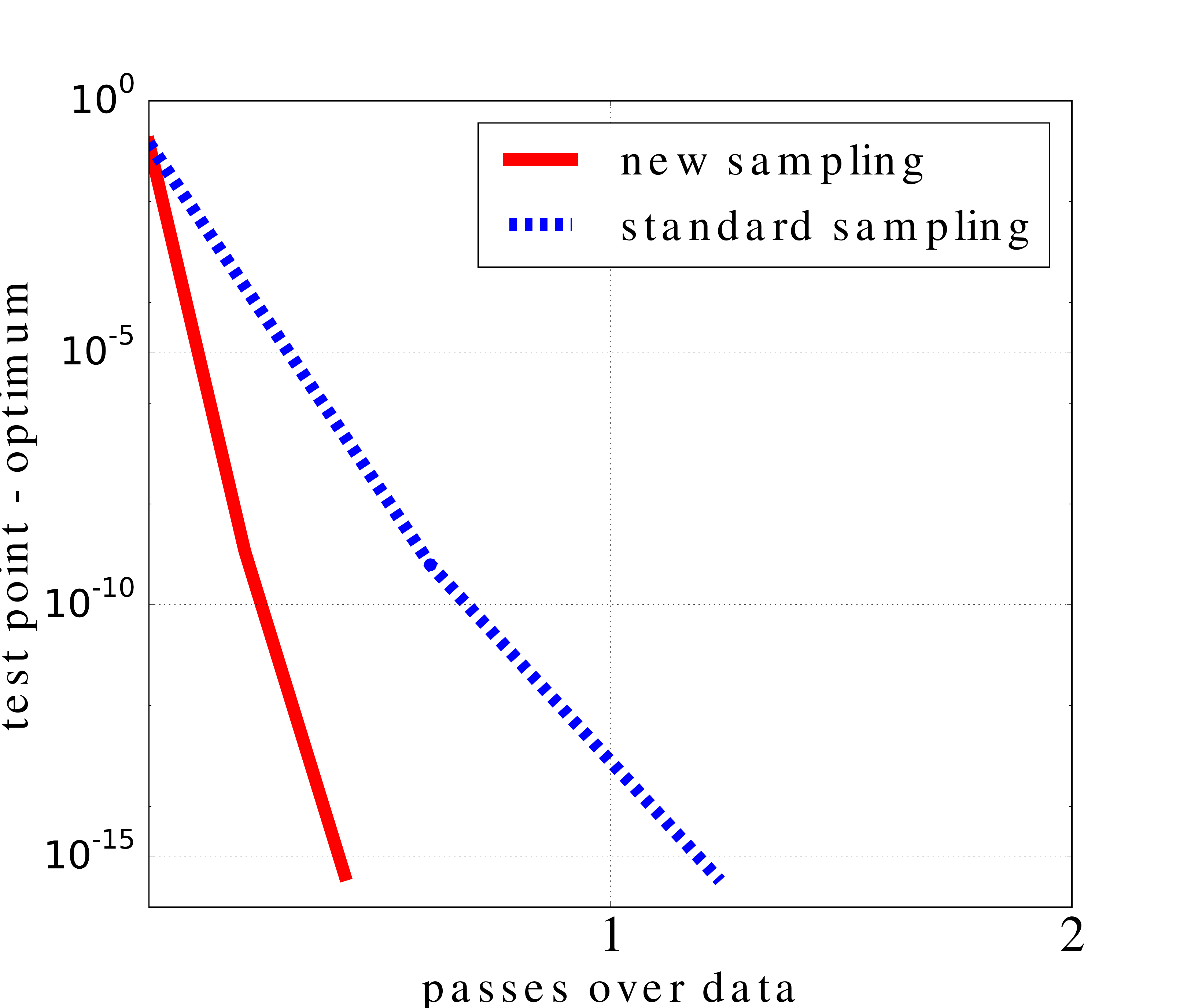}
            \caption[]%
            {{\small dorothea with $\tau=50$}}    
        \end{subfigure}
        \begin{subfigure}[b]{0.24\textwidth}  
            \centering 
            \includegraphics[width=\textwidth]{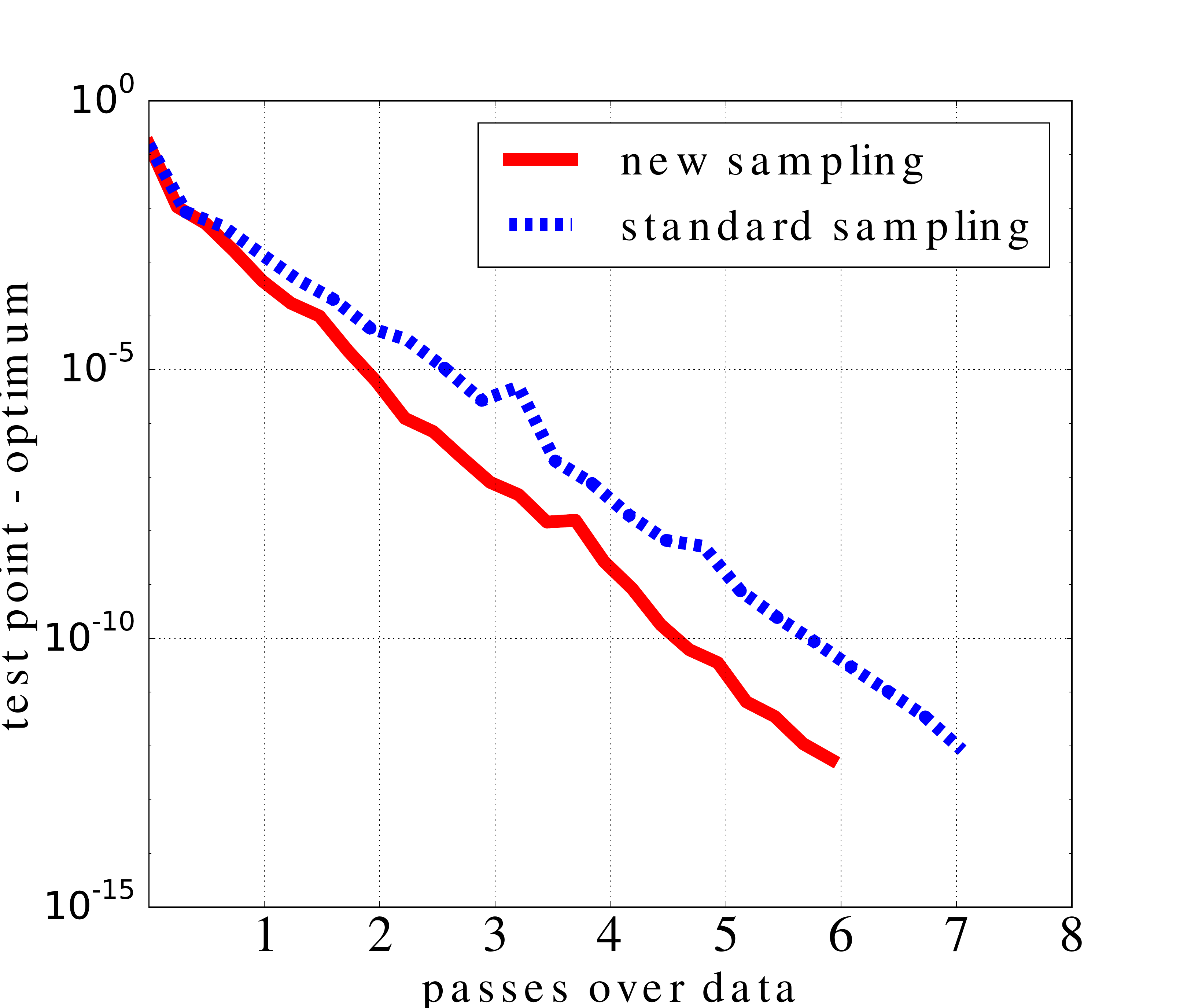}
            \caption[]%
            {{\small protein with $\tau=5$}}    
        \end{subfigure}
         \begin{subfigure}[b]{0.24\textwidth}  
            \centering 
            \includegraphics[width=\textwidth]{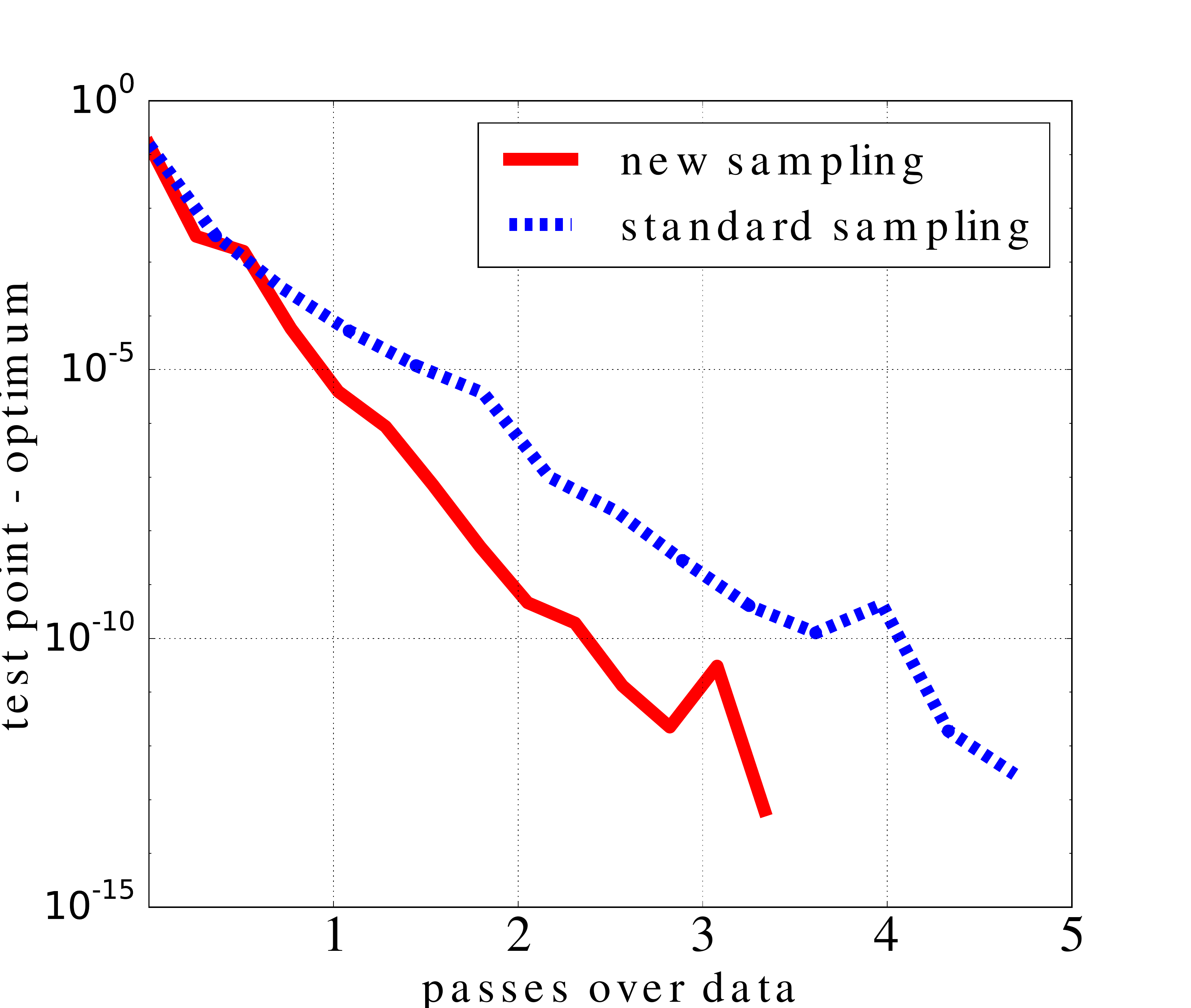}
            \caption[]%
            {{\small protein with $\tau=10$}}    
        \end{subfigure}
        \begin{subfigure}[b]{0.24\textwidth}  
            \centering 
            \includegraphics[width=\textwidth]{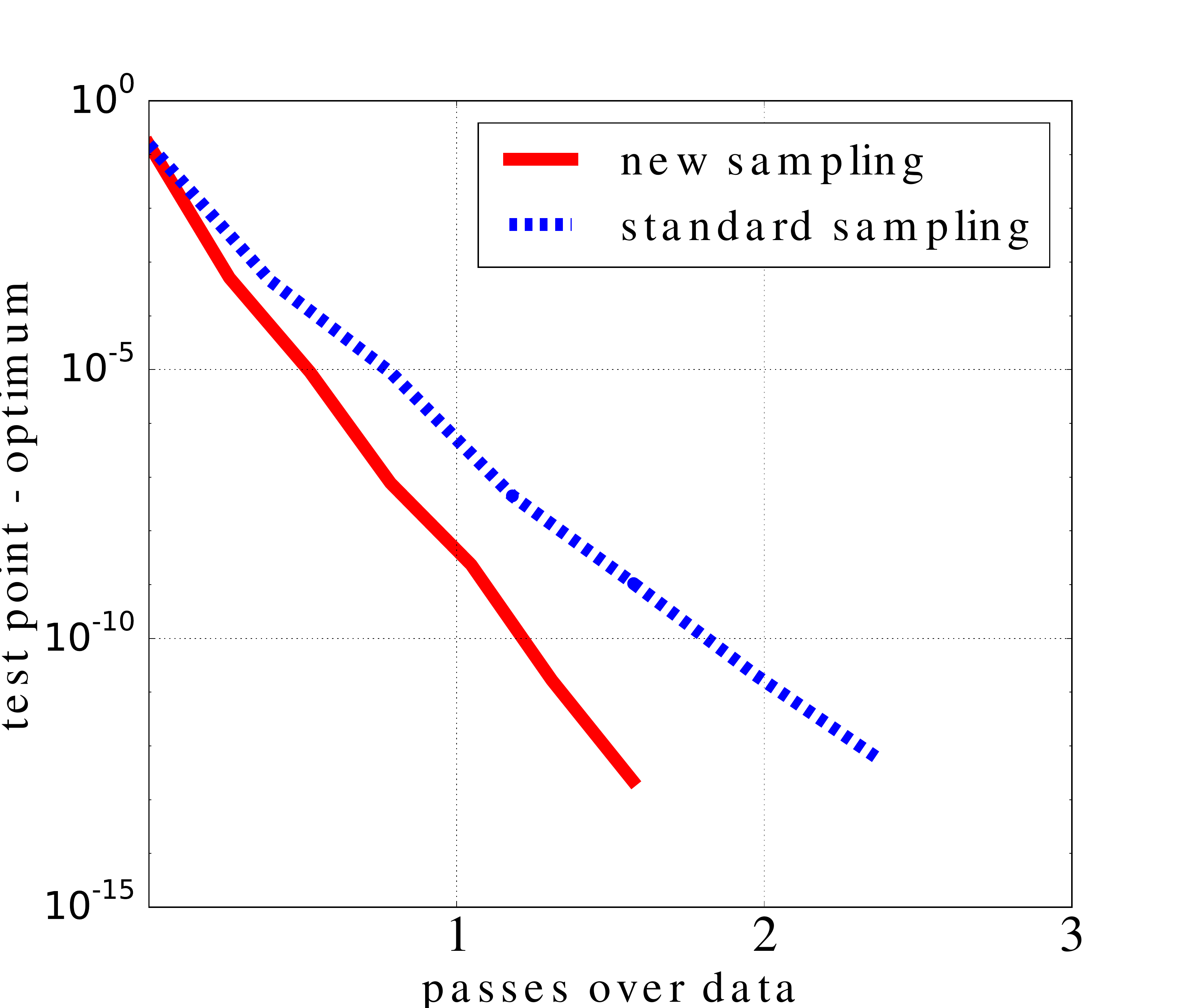}
            \caption[]%
            {{\small protein with $\tau=20$}}    
        \end{subfigure}
        \begin{subfigure}[b]{0.24\textwidth}  
            \centering 
            \includegraphics[width=\textwidth]{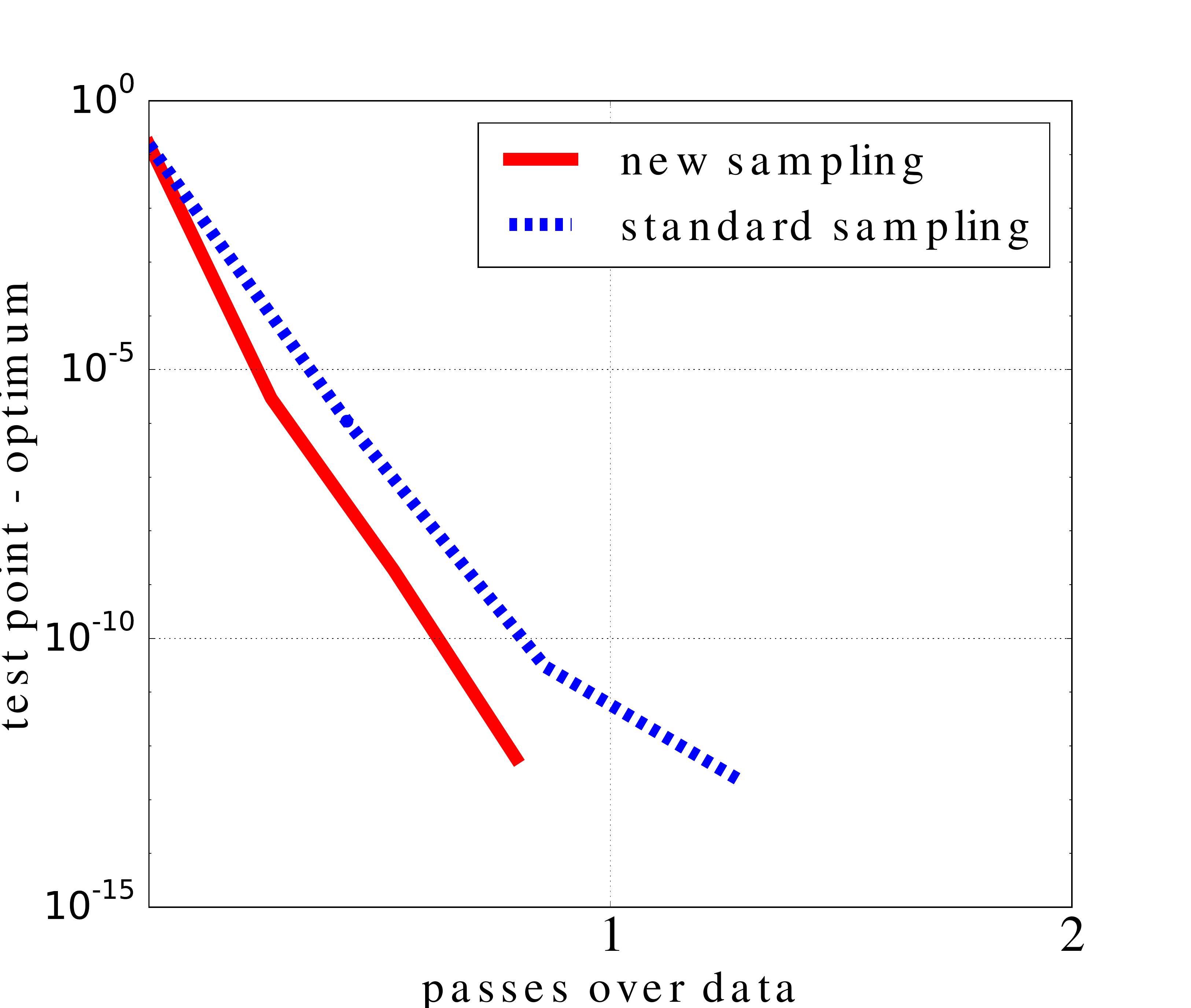}
            \caption[]%
            {{\small protein with $\tau=50$}}
            \label{fig:protein_50}    
        \end{subfigure}
        \caption[]
        {\small Logistic regression with $\lambda = 1/n$. Comparison between new and standard sampling with fine-tuned stepsizes for different values of $\tau$.}  
    \end{figure*}

\textbf{Experiment 2: New sampling vs. old sampling.} In Figure~\ref{fig:w8a_5} through Figure~\ref{fig:protein_50} we compare the performance of a standard parallel sampling against sampling of blocks $G_1, \dots, G_k$ output by Algorithm~\ref{alg:naive_chunks}. In each iteration we measure the time by \[\max_{i \in S_t} \{ \text{nnz}(A_i) \}\] and \[\max_{i \in [\tau]} \{ \text{nnz}(G_{(i)}) \}\] for the standard and new sampling respectively. This way we measure only the computations done by the core which is going to finish the last in each iteration, and consider the number of multiplications with nonzero entries of the data matrix as a proxy for time.

\section{Proofs}

As a first approximation, our proof is an extension of the proof of Shalev-Shwartz \cite{DualFreeSDCA} to accommodate an arbitrary sampling \cite{NSync,ALPHA,Quartz, SDNA}. For all $i$ and $t$ we let $u_i^{(t-1)} = -\nabla\phi_i (A_i^\top w^{(t)})$ and $z_i^{(t-1)} =  \alpha_i^{(t-1)} - u_i^{(t-1)}$.  We will use the following lemma.
\begin{lemma}[Evolution of $C_i^{(t)}$ and $B^{(t)}$] \label{lem:phiAiw_evol}
 For a fixed iteration $t$and all $i$ we have:
\begin{align}
\E_{\hat{S}}\left[ C_i^{(t-1)} - C_i^{(t)}\right] &= \theta \left[ \|\alpha_i^{(t-1)} - \alpha_i^* \|^2 - \|u_i^{(t-1)} - \alpha_i^*\|^2 + (1 - \theta p_i^{-1})\|z_i^{(t-1)}\|^2 \right] \label{eq:phiAiw_Ct}
\\ \E_{\hat{S}} \left[ B^{(t-1)} - B^{(t)} \right] &\geq  \frac{2\theta}{\lambda} (w^{(t-1)} - w^*)^\top \nabla P(w^{(t-1)}) - \frac{\theta^2}{ n^2 \lambda^2} \sum_{i=1}^n \frac{v_i}{p_i}\|z_i^{(t-1)}\|^2. \label{eq:phiAiw_Bt}
\end{align} 
\end{lemma}

{\em Proof.} It follows that for $i \in S_t$ using the definition \eqref{def:phiAitw_pot} we have
\begin{eqnarray*}
C_i^{(t-1)} - C_i^{(t)} &\overset{\eqref{def:phiAitw_pot}}{=}& \|\alpha_i^{(t-1)} - \alpha_i^*\|^2 - \|\alpha_i^{(t)} - \alpha_i^* \|^2
\\ &=& \|\alpha_i^{(t-1)} - \alpha_i^*\|^2 - \|(1-\theta p_i^{-1})(\alpha_i^{(t-1)} - \alpha_i^*) + \theta p_i^{-1} (u_i^{(t-1)} - \alpha_i^*)\|^2 
\\ &=& \| \alpha_i^{(t-1)} - \alpha_i^*\|^2 -(1-\theta p_i^{-1})\|\alpha_i^{(t-1)} - \alpha_i^*\|^2 - \theta p_i^{-1}\|u_i^{(t-1)} - \alpha_i^*\|^2 \\ &&+ \theta p_i^{-1}(1-\theta p_i^{-1})\|\alpha_i^{(t-1)} - u_i^{(t-1)}\|^2
\\ &=& \theta p_i^{-1}\left[ \|\alpha_i^{(t-1)} - \alpha_i^*\|^2 - \|u_i^{(t-1)} - \alpha_i^*\|^2 + (1-\theta p_i^{-1})\|z_i^{(t-1)}\|^2 \right]
\end{eqnarray*}
and for $i \notin S_t$ we have $C_i^{(t-1)} - C_i^{(t)} = 0$. Taking the expectation over $S_t$ we get the result.

For the second potential we get
\begin{eqnarray*}
B^{(t-1)} - B^{(t)} &\overset{\eqref{def:phiAitw_pot}}{=}& \|w^{(t-1)} - w^*\|^2 - \|w^{(t)} - w^*\|^2
\\ &=& \frac{2\theta}{n\lambda}\sum_{i\in S_t} p_i^{-1}(w^{(t-1)} - w^*)^\top A_iz_i^{(t-1)} - \frac{\theta^2}{n^2\lambda^2}\| \sum_{i \in S_t} p_i^{-1}A_i z_i^{(t-1)}\|^2.
\end{eqnarray*}
Taking the expectation over $S_t$, using inequality \eqref{eq:phiAiw_ESO}, and noting that 
\begin{equation} \label{eq:lemPw}
\frac{1}{n}\sum_{i=1}^n A_i z_i^{(t-1)} = \frac{1}{n}\sum_{i=1}^n A_i\nabla\phi (A_i^\top w^{(t-1)})+ \lambda w^{(t-1)} = \nabla P(w^{(t-1)}) ,
\end{equation} 
we get
\begin{align*}
\E \left[ B^{(t-1)} - B^{(t)} \right] &= \frac{2\theta}{n\lambda}\sum_{i=1}^n (w^{(t-1)} - w^*)^\top A_iz_i^{(t-1)} - \frac{\theta^2}{n^2\lambda^2}\E\left[ \|\sum_{i\in S_t}A_i (p_i^{-1}z_i^{(t-1)})\|^2 \right]
\\ & \stackrel{\eqref{eq:phiAiw_ESO}}{\geq}  \frac{2\theta}{n\lambda}\sum_{i=1}^n (w^{(t-1)} - w^*)^\top A_iz_i^{(t-1)} - \frac{\theta^2}{n^2\lambda^2}\sum_{i=1}^n p_i v_i \|p_i^{-1}z_i^{(t-1)}\|^2
\\ &\stackrel{\eqref{eq:lemPw}}{=} \frac{2\theta}{\lambda} (w^{(t-1)} - w^*)^\top \nabla P(w^{(t-1)}) - \frac{\theta^2}{n^2\lambda^2} \sum_{i=1}^n \frac{v_i}{p_i}\|z_i^{(t-1)}\|^2\qed
\end{align*}

\subsection{Proof of Theorem~\ref{thm:phiAiw_nonconvex} (nonconvex case)}

Combining \eqref{eq:phiAiw_Ct} and \eqref{eq:phiAiw_Bt}, we obtain
\begin{align*}
\E[D^{(t-1)} - D^{(t)}] &\geq \frac{\theta\lambda}{2n}\sum_{i=1}^n \frac{1}{L_i^2}\left[ \|\alpha_i^{(t-1)} - \alpha_i^* \|^2 - \|u_i^{(t-1)} - \alpha_i^*\|^2 + (1 - \theta p_i^{-1})\|z_i^{(t-1)}\|^2\right]
\\ &\qquad + \frac{\lambda}{2}\left[ \frac{2\theta}{\lambda} (w^{(t-1)} - w^*)^\top \nabla P(w^{(t-1)}) - \frac{\theta^2}{n^2\lambda^2}\sum_{i=1}^n \frac{v_i}{p_i}\|z_i^{(t-1)}\|^2 \right]
\\ &= \frac{\theta}{2n} \sum_{i=1}^n \left[ \frac{\lambda}{L_i^2} \left( C_i^{(t-1)} - \|u_i^{(t-1)} - \alpha_i^*\|^2\right) + \left( \frac{\lambda (1-\theta p_i^{-1})}{L_i^2} - \frac{\theta v_i}{n \lambda p_i}\right)\|z_i^{(t-1)}\|^2  \right]
\\ & \qquad + \theta (w^{(t-1)} - w^*)^\top \nabla P(w^{(t-1)})
\\ &\stackrel{\eqref{eq:phiAiw_eta}}{\geq} \frac{\theta}{2n} \sum_{i=1}^n \frac{\lambda}{L_i^2} \left( C_i^{(t-1)} - \|u_i^{(t-1)} - \alpha_i^*\|^2\right) + \theta (w^{(t-1)} - w^*)^\top \nabla P(w^{(t-1)}).
\end{align*}
Using \eqref{eq:phiAiw_smoothness} we have 
\[\|u_i^{(t-1)} - \alpha_i^* \|^2 = \|\nabla\phi_i(A_i^\top w^{(t-1)})~-~\nabla\phi_i(A_i^\top w^*)\|^2 \leq L_i^2 \|w^{(t-1)}~-~w^*\|^2.\] By strong convexity of $P$,  \[(w^{(t-1)} - w^*)^\top \nabla P(w^{(t-1)}) \geq  P(w^{(t-1)}) - P(w^*) + \frac{\lambda}{2}\|w^{(t-1)} - w^*\|^2\] and $P(w^{(t-1)}) - P(w^*) \geq \frac{\lambda}{2}\|w^{(t-1)} - w^*\|^2,$ which together yields \[(w^{(t-1)} - w^*)^\top \nabla P(w^{(t-1)}) \geq \lambda\| w^{(t-1)} - w^*\|^2.\] Therefore,
\begin{equation*}
\E[D^{(t-1)} - D^{(t)}] \geq \theta\left[ \frac{1}{n} \sum_{i=1}^n \frac{\lambda}{2L_i^2} C_i^{(t-1)} + \left(-\frac{\lambda}{2} + \lambda \right)B^{(t-1)} \right] = \theta D^{(t-1)}.
\end{equation*} 
It follows that
$
\E[D^{(t)}] \leq (1 - \theta)D^{(t-1)}
$, and repeating this recursively we end up with
$\E[D^{(t-1)}] \leq (1-\theta)^t D^{(0)} \leq e^{-\theta t}D^{(0)}. $
This concludes the proof of the first part of Theorem~\ref{thm:phiAiw_nonconvex}. The second part of the proof follows by observing that $P$ is $(L + \lambda)$-smooth, which gives $P(w) - P(w^*) \leq \frac{L + \lambda}{2}\|w - w^*\|^2$.

\subsection{Convex case}

For the next theorem we need an additional lemma:

\begin{lemma} \label{lem:convex_bound}
Assume that $\phi_i$ are $L_i$-smooth and convex. Then, for every $w$,
\begin{equation} \label{eq:phiAiw_convexlemma}
\frac{1}{n}\sum_{i=1}^n \frac{1}{L_i}\|\nabla\phi_i(w) - \nabla\phi_i(w^*)\|^2 \leq 2 \left( P(w) - P(w^*) - \frac{\lambda}{2}\|w - w^*\|^2 \right)
\end{equation} 
\end{lemma}
{\em Proof.}
Let $g_i(x) = \phi_i(x) - \phi_i(A_i^\top w^*) - \nabla\phi_i(A_i^\top w^*)^\top (x - A_i^\top w^*).$ Clearly, $g_i$ is also $l_i$-smooth. By convexity of $\phi_i$ we have $g_i(x) \geq 0$ for all $x$. It follows that $g_i$ satisfies $\|\nabla g_i (x)\|^2 \leq 2 l_i g_i(x).$ Using the definition of $g_i$, we obtain
\begin{align}
\|\nabla\phi_i (A_i^\top w) - \nabla\phi_i(A_i^\top w^*)\|^2 &= \|\nabla g_i(A_i^\top w)\|^2 \nonumber
\\ &\leq 2l_i [\phi_i(A_i^\top w)-\phi_i(A_i^\top w^*) - \nabla\phi_i (A_i^\top w^*)^\top (A_i^\top w - A_i^\top w^*)]. \label{lem:phiLi}
\end{align}
Summing these terms up weighted by $1/l_i$ and using \eqref{eq:phiAiw_optimality} we get
\begin{align*}
\frac{1}{n}\sum_{i=1}^n \frac{1}{l_i}\|\nabla\phi_i(A_i^\top w) - \nabla\phi_i(A_i^\top w^*)\|^2 
 &\stackrel{\eqref{lem:phiLi}}{\leq} \sum_{i=1}^n \frac{2}{n}[\phi_i(A_i^\top w) - \phi_i(A_i^\top w^*) - A_i\nabla\phi_i(A_i^\top w^*)^\top (w - w^*)]
\\ &\stackrel{\eqref{eq:phiAiw_optimality}}{=} 2\left[ P(w) - \frac{\lambda}{2}\|w\|^2 - P(w^*) + \frac{\lambda}{2}\|w^*\|^2 + \lambda {w^*}^\top (w - w^*) \right]
\\ &= 2 \left[ P(w) - P(w^*) - \frac{\lambda}{2}\|w - w^*\|^2 \right].\qed
\end{align*}

\subsection{Proof of Theorem~\ref{thm:phiAiw_convex}}

Combining \eqref{eq:phiAiw_Ct} and \eqref{eq:phiAiw_Bt}, we obtain
\begin{align*}
\E[E^{(t-1)} - E^{(t)}] &\geq  \frac{\theta}{n}\sum_{i=1}^n \frac{1}{2l_i}\left[ \|\alpha_i^{(t-1)} - \alpha_i^*\|^2 - \|u_i^{(t-1)} - \alpha_i^*\|^2 + (1-\theta p_i^{-1})\|z_i^{(t-1)}\|^2 \right]
\\ & \qquad + \frac{\lambda}{2}\left[\frac{2\theta}{\lambda} (w^{(t-1)} - w^*)^\top \nabla P(w^{(t-1)}) - \frac{\theta^2}{n^2\lambda^2}\sum_{i=1}^n \frac{v_i}{p_i}\|z_i^{(t-1)}\|^2\right]
\\ &= \frac{\theta}{n}\sum_{i=1}^n \left[ \frac{1}{2l_i} (C_i^{(t-1)} - \|u_i^{(t-1)} - \alpha_i^*\|^2) + \left( \frac{(1-\theta p_i^{-1})}{2l_i} - \frac{\theta v_i}{2 p_i \lambda n} \right) \right]
\\ &\qquad + \theta (w^{(t-1)} - w^*)^\top \nabla P(w^{(t-1)})
\\ &\stackrel{\eqref{eq:phiAiw_eta_convex}}{\geq} \frac{\theta}{n}\sum_{i=1}^n \left[ \frac{1}{2l_i} (C_i^{(t-1)} - \|u_i^{(t-1)} - \alpha_i^*\|^2) \right] + \theta (w^{(t-1)} - w^*)^\top \nabla P(w^{(t-1)})
\end{align*}

Using the convexity of $P$ we have $P(w^*) - P(w^{(t-1)}) \geq (w^{(t-1)} - w^*)^\top \nabla P(w^{(t-1)})$ and using Lemma~\ref{lem:convex_bound}, we have
\begin{eqnarray*}
\E[E^{(t-1)} - E^{(t)}] 
 &\stackrel{\eqref{eq:phiAiw_convexlemma}}{\geq} & \frac{\theta}{n} \sum_{i=1}^n \frac{1}{2l_i}C_i^{(t-1)} - \theta\left( P(w^{(t-1)}) - P(w^*) - \frac{\lambda}{2}\|w^{(t-1)} - w^*\|^2 \right) \\
 \qquad && + \theta (w^{(t-1)} - w^*)^\top \nabla P(w^{(t-1)})
\\ &\geq & \theta \left[ \frac{1}{n}\sum_{i=1}^n\frac{1}{2l_i}C_i^{(t-1)} + \frac{\lambda}{2}B^{(t-1)} \right]  \quad = \quad \theta E^{(t-1)}.
\end{eqnarray*}

This gives $\E[E^{(t)}] \leq (1 - \theta )E^{(t-1)}$, which concludes the first part of the Theorem~\ref{thm:phiAiw_convex}. The second part follows by observing, that $P$ is $(L + \lambda)$-smooth, which gives $P(w) - P(w^*) \leq \frac{L + \lambda}{2}\|w - w^*\|^2$.


\bibliography{DualFree} \bibliographystyle{unsrt}

\end{document}